%% file: ex_article.tex
\documentclass[hidelinks,onefignum,onetabnum]{siamart251104}

\input{ex_shared}

\ifpdf
  \hypersetup{
    pdftitle={Inhomogeneous Priors for Bayesian Inverse Problems},
    pdfauthor={Babak Maboudi Afkham, Tom\'as Soto, Mirza Karamehmedovi\'c and Lassi Roininen}
  }
\fi

\begin{document}

\maketitle

\begin{abstract}
    Many inverse problems arising in engineering and applied sciences involve unknown quantities with pronounced spatial inhomogeneity, such as localized defects or spatially varying material properties, making reliable uncertainty quantification particularly challenging. While Bayesian inverse problem methodologies provide a principled framework for assessing reconstruction reliability, commonly used Gaussian priors, such as Whittle–Mat\'ern models, impose globally homogeneous assumptions that limit their ability to capture such structure in large-scale settings. We introduce a new class of inhomogeneous priors defined via convolution with white noise, yielding nonstationary Whittle–Matérn–type random fields with a rigorous mathematical construction. These priors fit naturally within existing Bayesian well-posedness theory and enable efficient sampling by reducing prior realizations to the solution of a pseudo-differential equation, for which we develop numerical schemes with quantified approximation error. Numerical experiments in one-dimensional denoising and two-dimensional limited-angle X-ray tomography demonstrate significant improvements in reconstruction quality and uncertainty quantification, particularly in data-limited scenarios.
\end{abstract}

\begin{keywords}
Inverse Problems, Bayesian Framework, Priors, Pseudo-Differential Operators, Inhomogeneity
\end{keywords}

\begin{MSCcodes}
65J22, 62F15, 60G15, 60B11, 35S05
\end{MSCcodes}

\section{Introduction}

The most consequential scientific and engineering questions often rely on inferring hidden structure from indirect, incomplete, or noisy observations. Problems of this type, known as inverse problems, are central to medical imaging, where internal anatomy must be reconstructed from external measurements; to industrial monitoring, where early detection of defects can prevent catastrophic failure; and to scientific discovery, where physical parameters must be inferred from limited experimental data. As a result, the way inverse problems are modeled and solved has a direct impact on both the accuracy of reconstructions and the decisions based upon them.

Many inverse problems arising in applications involve unknown quantities that are inherently inhomogeneous. Examples include materials containing defects, biological tissues with spatially varying properties, and geological media exhibiting layered or fractured structures. In the renewable energy sector, wind turbine blades are susceptible to cracks, erosion, and surface corrosion induced by cyclic wind loading and environmental exposure; if left undetected, such defects can lead to measurable losses in energy output \cite{kong2023progress}. Similarly, in solar photovoltaic systems, corrosion of panel components and microstructural defects have been identified as primary challenges for long-term reliability \cite{dubchak2024modern}. In such problems, accurate inference requires not only estimating the unknown itself, but also resolving spatially varying structure and localized inhomogeneities that are critical for interpretation and decision-making.

Inverse problems are often ill-posed, in the sense that solutions may fail to exist or be unique, or may depend sensitively on fluctuations in the data. Regularization is therefore essential, effectively restricting the search space and encoding assumptions about the structure of the unknown. Classical regularization techniques typically impose global assumptions, such as uniform smoothness or sparsity, thereby implicitly treating the unknown as structurally homogeneous. While effective in many settings, these assumptions can be poorly matched to problems in which variability is localized or spatially varying. In such cases, localized heterogeneities may be smoothed out or underestimated, while noise may be incorrectly promoted to structural significance \cite{agapiou2024laplace,bhadra2021hallucinations}.

Bayesian formulations of inverse problems provide a principled framework for uncertainty quantification, in which regularization is encoded through prior distributions. From this perspective, priors translate modeling assumptions into probability distributions for the unknown, which directly affect both point estimates and uncertainty statements. While Bayesian methods enable the quantification of uncertainty associated with reconstructed features, this capability is fundamentally shaped by the choice of prior. In particular, most commonly used priors encode stationary or globally uniform assumptions on regularity or variability. As a consequence, such priors are unable to adapt to spatially varying structure and may oversmooth complex regions or under-regularize simpler ones, thereby distorting both reconstructions and associated uncertainty assessments.

A natural way to formulate Bayesian inverse problems in a manner that is robust to discretization is to work directly at the level of function spaces, rather than through finite-dimensional parameterizations. In this setting, the unknown is modeled as a random function on a continuum, and prior and posterior distributions are defined independently of any particular discretization. This infinite-dimensional perspective provides a mathematically consistent framework for uncertainty quantification and ensures stability of inference under mesh refinement. Importantly, while infinite-dimensional formulations enable spatially resolved modeling in principle, the ability to capture spatially varying structure in practice depends critically on the expressiveness of the chosen prior distribution.

With the growing importance of uncertainty quantification, there is increasing demand for solving high-dimensional inverse problems within this infinite-dimensional Bayesian framework. However, standard sampling methods for exploring posterior distributions, such as Markov chain Monte Carlo (MCMC) or Metropolis–Hastings, typically scale poorly with dimension. The theory of infinite-dimensional Bayesian inverse problems, pioneered in \cite{Stuart_2010,Dashti2017,Lasanen-2012-II}, addresses this challenge by enabling the design of discretization-invariant MCMC algorithms whose performance remains stable under mesh refinement, thereby facilitating efficient exploration of high-dimensional posterior distributions \cite{cotter2013mcmc,boehl2022ensemble}.

Despite the success of infinite-dimensional Bayesian formulations, a central challenge remains the construction of suitable prior distributions on function spaces. A widely used approach is to specify a covariance kernel, which induces a covariance operator for a Gaussian prior \cite{seeger2004gaussian,rue2005gaussian,mandrekar2015stochastic,ibragimov2012gaussian}. In large-scale settings, however, explicit representations of these operators can be computationally prohibitive. The Karhunen–Lo\`eve (KL) expansion provides a scalable alternative, representing random functions through their spectral decomposition and allowing priors to be defined implicitly without forming full covariance operators.

Whittle–Matérn priors \cite{matern1986spatial,whittle1963stochastic,lindgren11spde,roininen2014whittle,lindgren21gaussian,porcu2024matern} form a broad and expressive class of Gaussian priors whose covariance structure and spectral properties are well understood. Governed by a correlation length, a smoothness exponent, and a variance parameter, these priors provide flexible control over global regularity and scale. However, because these parameters are spatially uniform, Whittle–Matérn priors encode globally homogeneous and isotropic assumptions, causing the resulting random fields to exhibit uniform statistical properties across the domain. Capturing spatially varying inhomogeneity or anisotropy therefore requires extensions beyond classical formulations and remains an active area of research.

In this paper, we introduce a class of inhomogeneous priors for Bayesian inverse problems, whose distributions are defined via convolution with white noise. Within this framework, spatially varying Whittle–Matérn covariance operators arise naturally, providing a rigorous construction of nonstationary Gaussian random fields. Moreover, this formulation is directly compatible with existing Bayesian well-po\-sed\-ness theory \cite{Stuart_2010,latz2020well,Lasanen-2012-II}, allowing previously studied inverse problems to be incorporated without the need for further theoretical redevelopment.

We show that sampling realizations of the proposed random variables reduces to solving a stochastic partial differential equation (SPDE) characterized by a pseudo-differential operator \cite{hormander1985analysis3,grubb2012functional}. To this end, we develop a detailed numerical scheme and quantify its approximation error, an aspect that has received comparatively limited attention in the existing literature. Moreover, when the symbol of the pseudo-differential operator exhibits rapid decay, as observed in all numerical examples presented in this paper, the sampling procedure simplifies considerably and can be interpreted as a straightforward modification of existing algorithms for Whittle–Matérn priors. We provide a complete and implementable algorithm for sampling from these inhomogeneous priors.

Designing solvers and MCMC methods for large-scale inverse problems remains a significant challenge. We show that the proposed inhomogeneous priors can be naturally integrated with modern and efficient inference algorithms, including the No-U-Turn Sampler (NUTS) for Hamiltonian Monte Carlo (HMC) and the Limited-memory Broyden–Fletcher–Goldfarb–Shanno (L-BFGS) method for computing maximum a posteriori (MAP) estimates.

We assess the performance of our approach on a one-dimensional denoising problem and a two-dimensional limited-angle X-ray computed tomography example. These experiments demonstrate that spatially varying priors are particularly advantageous in data-limited settings, where adaptivity to localized structure is essential. The numerical results indicate substantial improvements in reconstruction accuracy compared with standard methods, while simultaneously providing reliable uncertainty quantification.

The paper is organized as follows. In \Cref{sec:bayesian-inverse-problems} we review Bayesian inverse problems with priors defined on a continuum. In \Cref{ss:spde-pdos} we introduce a class of spatially inhomogeneous random functions as solutions to a particular SPDE and formalize their laws as prior distributions for Bayesian inverse problems. In \Cref{sec:computational-method} we present a computational algorithm for sampling from these inhomogeneous priors. In \Cref{eq:numerical-results} we report numerical results for a one-dimensional denoising problem and a limited-angle two-dimensional X-ray computed tomography problem, and investigate the performance of the proposed method in exploring posterior distributions constructed with inhomogeneous priors. Concluding remarks are given in \Cref{sec:conclusion}.

\section{Bayesian Inverse Problems on a Continuum}\label{sec:bayesian-inverse-problems}
We let normal letters $x$ stand for scalars, bold letters $\boldsymbol x$ for vectors and capital letters $X$ for random variables. The spaces in which random variables take realizations are explicitly determined, e.g.~if $X$ takes values in a space $\mathcal{X}$, then we say that $X$ is an $\mathcal{X}$-valued random variable.

Let $H$ be a Hilbert space. In this paper we consider inverse problems of the type
\begin{equation}\label{eq:inverse}
    \boldsymbol{y} = \mathcal{G}(\xi) + \boldsymbol{e},
\end{equation}
where $\xi\in H$ is the unknown in the inverse problem, $\mathcal{G} : H \to \mathbb{R}^m$ is the forward operator, $\boldsymbol e \in \mathbb{R}^m$ is an additive noise vector, and $\boldsymbol y\in \mathbb{R}^m$ is the measurement vector. The inverse problem comprises inferring the unknown $\xi$ in \eqref{eq:inverse} from the noisy measurement $\boldsymbol y$.

In the Bayesian formulation of inverse problems, we recast \eqref{eq:inverse} in a statistical framework by modeling the unknown, noise and measurements as random variables. Let $X$ be an $H$-valued random variable presenting the unknown $\xi$, and $Y$ and $E$ be $\mathbb{R}^m$-valued random variables, corresponding to $\boldsymbol y$ and $\boldsymbol e$, respectively. The statistical analogue of \eqref{eq:inverse} then takes the form
\begin{equation} \label{eq:statistical_inverse}
    Y = \mathcal{G}(X) + E.
\end{equation}
Here $\mathcal{G}$ is treated as a deterministic map, often representing the physical or mathematical model underlying the measurement. The goal of the statistical inverse problem \eqref{eq:statistical_inverse} is to characterize the conditional distribution of $X$ given observations $Y=\, \boldsymbol y$, written as $X|Y=\boldsymbol y$. This distribution is referred to as the \emph{posterior}. 

Bayes’ theorem \cite{kaipio2005statistical} traditionally provides the standard framework for decomposing the posterior distribution into more tractable components. However, since $X$ in our setting is function-valued, one generally cannot represent the law of $X$ or of $X|Y$ by a density with respect to a canonical dominating measure (e.g.~the Lebesgue measure on an Euclidean space). Instead, these distributions are formulated in terms of probability measures on infinite-dimensional function spaces.

Let $\mu_X$ and $\mu_{X|Y=\, \boldsymbol y}$ denote the prior and the posterior probability measures, respectively. The measure-theoretic formulation of Bayes' theorem states that
\begin{equation}
    \frac{\mathrm{d} \mu_{X|Y=\, \boldsymbol y}}{\mathrm{d} \mu_X}(\xi) = \frac{1}{c} L(\boldsymbol y|\xi),
\end{equation}
where $L(\boldsymbol y|\xi)$ is the likelihood of observing the measurement $\boldsymbol y$ given the unknown $\xi$, and $c\in \mathbb{R}^+$ is the normalization constant. Here $\mathrm{d}{\cdot} / \mathrm{d}{\cdot}$ is the \emph{Radon-Nikodym derivative} \cite{ibragimov2012gaussian} between measures.
Details and references for this formulation of the posterior distribution are given in \Cref{ss:functional-bayes-ip} below.

In the next section, we summarize common approaches for constructing probability measures on function spaces, which serve as prior distributions. We then discuss how realizations of such measures can be interpreted as solutions to certain fractional SPDE's. To conclude this section, we present a compact formulation of the statistical inverse problem, listing the key random variables and their dependencies:
\begin{equation} \label{eq:inverse-summary}
    \begin{aligned}
        Y &= \mathcal{G}(X) + E, \\
        X &\sim \mu_X, \\
        E &\sim \mathcal{N}(0,\Sigma_{\text{noise}}),\\
        Y|( X=\xi ) &\sim \mathcal{N}(\mathcal{G}(\xi),\Sigma_{\text{noise}} ).
    \end{aligned}
\end{equation}
This summary will become increasingly relevant as we introduce additional complexity and hierarchical structure into the inverse problem.

\subsection{Gaussian Random Functions}
Consider the probability space $(H, \mathcal{B(H)},\allowbreak \mathbb{P})$, with $H$ a Hilbert space, $\mathcal B(H)$ the Borel $\sigma$-algebra of $H$ and $\mathbb P$ a Borel probability measure. The prior distributions for $X$, introduced in the previous section, will be of this form. We assume here that $H=L^2(\Omega)$, where $\Omega \subset \mathbb{R}^d$ is a domain, equipped with the standard $L^2$ inner product $\langle \cdot,\cdot \rangle_H$ and norm $\| \cdot\|_H$.

We say that a Borel-measurable random element $X$ on $H$ is an $H$-valued Gaussian random function on $\Omega$ if, for any $h \in H$, the $\mathbb{R}$-valued random variable $\langle X,h\rangle_H$ is Gaussian. Equivalently, its characteristic function takes the form \cite{ibragimov2012gaussian}
\begin{equation}\label{eq:gaussian-hilbertspace}
    \phi_X(h) \mathrel{\mathop:}= \mathbb{E} \left[ \exp (\boldsymbol{i}\langle X, h \rangle_H) \right] = \exp\left(\boldsymbol{i}\langle m, h \rangle_H -\frac{1}{2} \langle \mathcal C h,h\rangle_H \right),
\end{equation}
where $\boldsymbol{i} \in \mathbb{C}$ is the imaginary unit, $m\in H$ is a mean function and $\mathcal C: H\to H$ is a covariance operator that is bounded, self-adjoint, positive and trace-class. In this case, we denote the law of $X$ on $H$ by $\mu_X \mathrel{\mathop:}= \mathbb P \circ X^{-1} \mathrel{\mathop:}= \mathcal N(m, \mathcal C)$, and write $X \sim \mathcal \mu_X$.

The following result fully characterizes Gaussian random functions in terms of the spectral decomposition of $\mathcal C$.

\begin{proposition} \label{thm:kl}
    \cite{kuo2006gaussian,ibragimov2012gaussian} Let $m\in H$ and $\mathcal C:H\to H$ be a mean function and a covariance operator as stated above. Let $\{e_j\}_{j=1}^{\infty}$ be the eigenfunctions and $\{ \lambda_j \}_{j=1}^{\infty}$ be the corresponding eigenvalues of $\mathcal C$, sorted in decreasing order. Then, $X\sim \mathcal N(m,\mathcal C)$ if and only if $X$ has the series expansion
    \begin{equation} \label{eq:kl}
        X = m + \sum_{j=1}^{\infty} \sqrt{\lambda_j} X_j e_j,
    \end{equation}
    where $X_j\sim \mathcal N(0,1)$ are independent real-valued random variables for all $j \geq 1$. The infinite series is interpreted in the sense that $(\sqrt{ \lambda_j })_{j \geq 1} \in \ell^2$, and so $\mathbb{E} \left[\| X -m \|^2_{H} \right] < \infty$.
\end{proposition}

The expansion \eqref{eq:kl} is known as the \emph{Karhunen-Lo\`eve} (KL) expansion \cite{ibragimov2012gaussian} of $X$. \Cref{thm:kl} means that, in order to define a Gaussian random function as a prior, one may either construct a covariance operator $\mathcal C$ directly or specify its spectral decomposition. In practice, $X$ is often approximated by truncating the KL expansion, yielding $\widetilde X$, such that the approximation error satisfy $\mathbb{E} \left[\| X - \widetilde X \|^2_{H} \right] < \delta_{\text{KL}}$, i.e.~variance loss of at most $\delta_{\text{KL}}$ is tolerated.

Both of these tasks, designing $\mathcal C$ or its spectral decomposition, are challenging in high-dimensional inverse problems. Discretizations of a covariance operator often produce large covariance matrices and exploring a posterior of the form \eqref{eq:inverse-summary} often requires costly matrix inversion and factorization. On the other hand, working directly with the spectral decomposition also requires specifying infinitely many basis functions which is difficult to implement in a consistent and computationally tractable manner.

A promising approach \cite{lindgren11spde, lindgren21gaussian} that can alleviate these difficulties is viewing $X$ as a solution of an elliptic stochastic SPDE. To illustrate this, let $\Omega = [0,1)$ be a periodic domain (the $1$-dimensional torus) and $k \geq 1$ an integer. Set $m\equiv0$, and define the eigenvalues and eigenfunctions as $\lambda_j \mathrel{\mathop:}= 1/(\boldsymbol{i}2\pi j)^{2k}$ for $j\in \mathbb N^+$ and $e_j(x) \mathrel{\mathop:}= \exp(\boldsymbol{i} 2\pi j x)$ (the Fourier basis with integer frequencies).

In this case, the KL expansion \eqref{eq:kl} becomes
\begin{equation}
    X = \sum_{j=1}^{\infty} \frac{1}{(\boldsymbol{i}2\pi j)^k} X_j \exp( \boldsymbol{i} 2\pi j x ).
\end{equation}
Projection $X$ onto $e_j(x)$, i.e.~setting $\widehat{X}(j) \mathrel{\mathop:}= \langle X, e_j\rangle_H$, gives
\[
    \widehat{X}(j) = \frac{1}{(\boldsymbol{i}2\pi j)^k} X_j \qquad \implies \qquad (\boldsymbol{i}2\pi j)^k  \widehat{X}(j) = X_j.
\]
As multiplication by $-(\boldsymbol{i}2\pi j)^k$ in the Fourier domain corresponds to taking $k$th derivative in the physical domain, this relation is equivalent to the differential equation
\begin{equation} \label{eq:SPDE_simple}
    -\frac{\partial^k }{\partial x^k} X = \Psi,
\end{equation}
where the source term $\Psi$ is a tempered distribution defined by its Fourier transform, $\widehat{\Psi}(j) \mathrel{\mathop:}= X_j$, and is referred to as \emph{Hilbert-space white noise}. See \Cref{ap:whitenoise} for details.

In summary, the solution of the SPDE \eqref{eq:SPDE_simple} is a Gaussian random function that admits a KL-expansion of the form \eqref{eq:kl}. The associated covariance operator, as characterized by \Cref{thm:kl}, takes the form $\mathcal C = ( \partial^k/\partial x^k )^{-2}$. This further highlights the complexity of covariance operators for such Gaussian random functions.

Designing elliptic PDEs and developing efficient solvers has been a central theme of successful research over the past century. This body of work provides a powerful foundation for constructing new types of priors in inverse problems. In particular, it is well known \cite{lindgren11spde, roininen2014whittle} that Whittle–Mat\'ern covariance operators \cite{genton2001classes} can be reformulated as solutions to SPDE's of the form
\begin{equation} \label{eq:whittle_matern}
    c(\sigma,\alpha)(\sigma - \Delta)^{\alpha/2} X = \Psi,
\end{equation}
ensuring $X$ is a Gaussian random function with KL expansion as in \Cref{thm:kl}. This formulation offers practical flexibility for designing rich classes of Gaussian prior models. To the best of the authors’ knowledge, however, all priors of this type in previous literature are homogeneous, i.e.~$c, \sigma, \alpha$ and $\Delta$ are independent of the spatial variable $\boldsymbol x$.

In the following sections, we introduce a new class of inhomogeneous covariance operators that can serve as priors for inverse problems. Specifically, we generalize \eqref{eq:whittle_matern} via a pseudo-differential operator formulation, and we show that the resulting solution $X$ defines a well-posed probability measure $\mu_X$. We then employ this prior in an inverse problem setting and demonstrate its enhanced flexibility.

\subsection{Spatially Inhomogeneous SPDEs}\label{ss:spde-pdos}

Here we introduce an extension of the SPDE \eqref{eq:whittle_matern} with a \emph{spatially inhomogeneous} length scale function $\sigma$. Since we will extensively work with functions and operators in terms of their Fourier transforms, we fix $\mathbb{T}^d \mathrel{\mathop:}= [0, 1)^d$, the $d$-dimensional torus, as the domain of the SPDE. 

We refer to \Cref{ap:whitenoise} for an explicit construction of the toroidal white noise $\Psi$, and a brief overview of its functional-analytic properties. Crucially, $\Psi$ is a random tempered distribution, and a Gaussian field in the sense that the dual pairing
\[
    \langle \Psi, f \rangle_{\mathcal{D}'(\mathbb{T}^d) \times C^{\infty}(\mathbb{T}^d)}
\]
is a centered normally distributed random variable with variance $\|f\|_{L^2(\mathbb{T})}^2$ for every real-valued test function $f \in C^{\infty}(\mathbb{T}^d)$. For any arbitrarily small $\epsilon > 0$, a version of $\Psi$ exists as a random element in the negative-order Sobolev space $H^{-d/2-\epsilon}(\mathbb{T}^d)$  (see \Cref{ss:whitenoise-sobolev} for details).

We then consider the SPDE
\begin{equation}\label{eq:SPDE_inhomogeneous}
    p_{\alpha}(\cdot, D) X = \Psi
\end{equation}
for some $\alpha > 0$, where $p_{\alpha}(\cdot, D)$ stands for the \emph{pseudo-differential operator} given by the symbol $p_{\alpha}(\cdot, \cdot) \colon \mathbb{T}^d \times \mathbb{Z}^d \to \mathbb{C}$, defined as
\begin{equation}\label{eq:PDO-symbol}
    \mathbb{T}^d \times \mathbb{Z}^d \owns (x, \eta) \mapsto c(\sigma, \alpha) \left( \sigma(x) + |\eta|^2 \right)^{\alpha/2}.
\end{equation}
More precisely, under suitable assumptions on the length scale function $\sigma$ (which we will discuss below), the action of $p_{\alpha}(\cdot, D)$ on a test function $f \in C^{\infty}(\mathbb{T}^d)$ is given by
\begin{equation}\label{eq:pdo-definition}
    p_{\alpha}(\cdot, D) f (x) \mathrel{\mathop:}= \sum_{\eta \in \mathbb{Z}^d} p_{\alpha}(x, \eta) \widehat{f}(\eta) e^{ \boldsymbol{i} 2\pi \eta \cdot x },
\end{equation}
and $p_{\alpha}(\cdot, D)$ extends by duality as a continuous operator on the space $\mathcal{D}'(\mathbb{T}^d)$ of tempered distributions (see \Cref{ss:whitenoise-pdo} for details). We refer to e.g.~\cite[Chapters 3 and 4]{Ruzhansky_Turunen} more information on periodic analysis and pseudo-differential operators on $\mathbb{T}^d$.

We note that, in accordance with the discussion surrounding \eqref{eq:SPDE_simple}--\eqref{eq:whittle_matern}, the operator $p_{\alpha}(\cdot, D)$ defined in terms of the symbol \eqref{eq:PDO-symbol} can be thought of as a variant of the differential operator on the left-hand side of \eqref{eq:whittle_matern} with a spatially varying length scale. The SPDE \eqref{eq:SPDE_inhomogeneous} is in this sense a spatially inhomogeneous generalization of \eqref{eq:whittle_matern}. We also note that although the operation of a pseudo-differential operator is defined in terms of the Fourier transform of a complex-valued test function, a symbol like \eqref{eq:PDO-symbol} \textbf{guarantees} that $p(\cdot, D)$ maps real-valued test functions (resp.~distributions) to real-valued test functions (resp.~distributions); see \Cref{ss:whitenoise-pdo}.

Before discussing the solutions of \eqref{eq:SPDE_inhomogeneous}, let us clarify the assumptions on the length scale function $\sigma$. In order to appeal to general results concerning pseudo-differential operators, we will work exclusively with symbols in the H\"or\-man\-der classes $S^m( \mathbb{T}^d \times \mathbb{Z}^d )$, $m \in \mathbb{R}$. For our purposes, it suffices to know \cite[Theorem 4.5.3]{Ruzhansky_Turunen} that a toroidal symbol $p \colon \mathbb{T}^d \times \mathbb{Z}^d \to \mathbb{C}$ is in the class $S^m( \mathbb{T}^d \times \mathbb{Z}^d )$ if and only if $p$ can be extended as a smooth function defined in $\mathbb{T}^d \times \mathbb{R}^d$ so that
\[
    \sup_{(x, \eta) \in \mathbb{T}^d \times \mathbb{R}^d} (1 + |\eta| )^{|\beta| - m} \left| \partial^{( \beta )}_{\eta} \partial^{( \gamma )}_{x} p(x, \eta) \right| < \infty
\]
for all multi-indices $\beta$, $\gamma \in \mathbb{N}_0^d$. We can thus see that the symbol $p_{\alpha}(\cdot, \cdot)$ given by \eqref{eq:PDO-symbol} is in $S^{\alpha}( \mathbb{T}^d \times \mathbb{Z}^d )$ if $\sigma \colon \mathbb{T}^d \to (0, \infty)$ is \emph{smooth and uniformly bounded from above and below}, which we will assume in the sequel.

Thus, because pseudo-differential operators with symbol in $S^m( \mathbb{T}^d \times \mathbb{Z}^d )$ are bounded from $H^s(\mathbb{T}^d)$ to $H^{s-m}(\mathbb{T}^d)$ for all $m$, $s \in \mathbb{R}$ \cite[Corollary 4.8.3]{Ruzhansky_Turunen}, we can formally invert the pseudo-differential operator on the left-hand side of \eqref{eq:SPDE_inhomogeneous} to get a solution
\begin{equation} \label{eq:SPDE_solution}
    X = \left[ p_{\alpha}(\cdot, D) \right]^{-1} \Psi,
\end{equation}
the right-hand side being a random element in the Sobolev space $H^{\alpha-d/2-\epsilon}(\mathbb{T}^d)$. If $\alpha > d$, the Sobolev embedding theorem would further imply that the random field $X$ has H\"older $r$-regular sample functions for any $r \in (0, \alpha - d)$; see \Cref{ss:whitenoise-pdo}.

However, the inversion of spatially inhomogeneous pseudo-differential operators is in general a very difficult question and well beyond the scope of this work, so we will instead work with approximate inverses or \emph{parametrices} of the operators $p_{\alpha}(\cdot, D)$.

To this end, we note that for every $\alpha > 0$ and admissible length scale function $\sigma$, the symbol given by \eqref{eq:PDO-symbol} is \emph{elliptic} in that, besides being in the class $S^{\alpha}( \mathbb{T}^d \times \mathbb{Z}^d )$, it also satisfies
\[
    \inf_{(x, \eta) \in \mathbb{T}^d \times \mathbb{Z}^d} (1 + |\eta|)^{-\alpha} | p_{\alpha}(x, \eta) | > 0.
\]
This means \cite[Theorem 4.9.6]{Ruzhansky_Turunen} that, for any $N \in \mathbb{N}_0 \cup \{\infty\}$, there exists an operator $q_{-\alpha}^{(N)}(\cdot, D)$ with symbol
\[
    q_{-\alpha}^{(N)}( \cdot, \cdot ) \mathrel{\mathop:}= \sum_{j=0}^N q_{-\alpha}^j(\cdot, \cdot) \in S^{-\alpha}( \mathbb{T}^d \times \mathbb{Z}^d ),
\]
such that
\[
    q_{-\alpha}^{(N)}( \cdot, D ) \circ p_{\alpha}( \cdot, D ) - I,
\]
is for every $N$ a pseudo-differential operator with symbol in $S^{-N-1}( \mathbb{T}^d \times \mathbb{Z}^d )$. More precisely, we can take $q_{-\alpha}^{0}(\cdot, \cdot) = p_{\alpha}(\cdot, \cdot)^{-1}$, and we refer to \cite[Theorem 4.9.13]{Ruzhansky_Turunen} for a recursive formula for the higher-order terms $q_{-\alpha}^{j}(\cdot, \cdot)$, $j \geq 1$.

So, if $X$ is a solution to the SPDE \eqref{eq:SPDE_inhomogeneous} (which by ellipticity exists at least modulo an infinitely smooth term as a random element in $H^{\alpha-d/2-\epsilon}(\mathbb{T}^d)$), then we can approximate $X$ modulo a smoothing term by $q_{-\alpha}^{(N)}(\cdot, D) \Psi$, in the sense that
\begin{equation} \label{eq:approx-error}
    X - q_{-\alpha}^{(N)} (\cdot, D) \Psi \in H^{\alpha-d/2-\epsilon+N+1}(\mathbb{T}^d)
\end{equation}
with full probability. This can already for $N = 0$ and $\alpha > d$ relatively close to $d$ be very useful in reconstructions of rough features, as we will see in \Cref{sec:de-noising,sec:hierarchical}.

In summary, we have the following structural properties for the approximate solutions to the inhomogeneous SPDE \eqref{eq:SPDE_inhomogeneous}. Details concerning the functional-analytic properties and the Gaussian covariance structure mentioned here can be found in \Cref{ss:whitenoise-sobolev,ss:whitenoise-pdo}.

\begin{proposition}
    Let $q_{-\alpha}(\cdot, D) \mathrel{\mathop:}= q_{-\alpha}^{(N)}(\cdot, D)$, $N \in \mathbb{N}_{0} \cup \{\infty\}$, be any of the approximate parametrices for the pseudo-differential operator $p_{\alpha}(\cdot, D)$ in \eqref{eq:SPDE_inhomogeneous} discussed above, with corresponding symbol $q_{-\alpha}(\cdot, \cdot)$. Assume that $\alpha > d$.

    Then the approximate solution $q_{-\alpha}(\cdot, D) \Psi$ to the SPDE \eqref{eq:SPDE_inhomogeneous} exists as a random element in the H\"older space $C^{r}(\mathbb{T}^d)$ for any $r \in (0, \alpha-d)$, and $q_{-\alpha}(\cdot, D) \Psi$ is a centered Gaussian field with covariance operator
    \[
        \mathcal{C} \mathrel{\mathop:}= q_{-\alpha}(\cdot, D) \circ q_{-\alpha}^{\intercal}(\cdot, D)
    \]
    in the sense of \eqref{eq:gaussian-hilbertspace}, where $q_{-\alpha}^{\intercal}(\cdot, D)$ stands for the transpose of the operator $q_{-\alpha}(\cdot, D)$ with respect to the dual pairing between $\mathcal{D}'(\mathbb{T}^d)_{|C^{\infty}(\mathbb{T}^d)}$ and $C^{\infty}(\mathbb{T}^d)$.

    We have the KL-like expansion
    \[
        q_{-\alpha}(\cdot, D)\Psi(x) = \sum_{ \eta \in \mathbb{Z}^d } q_{-\alpha}(x, \eta) w_{\eta} e^{\boldsymbol{i} 2\pi \eta \cdot x},
    \]
    where the $w_{\eta}$ are the Fourier coefficients of the white noise $\Psi$, with almost sure pointwise convergence, and with the error estimate
    \[
        \sup_{x \in \mathbb{T}^d} \mathbb{E} \left[ \left| q_{-\alpha}(\cdot, D)\Psi(x) - \sum_{ |\eta| \leq M } q_{-\alpha}(x, \eta) w_{\eta} e^{\boldsymbol{i} 2\pi \eta \cdot x} \right|^2 \right]
        \leq C_{q, \alpha, d} M^{d - 2\alpha}.
    \]
\end{proposition}
The constant $C_{q, \alpha, d}$ appearing above can be reasonably estimated in the setting of our applications; see \Cref{ss:whitenoise-pdo}.
\subsection{Setup of the Functional Bayesian Inverse Problem}\label{ss:functional-bayes-ip}

As a quick recap of the introduction, recall that we consider the statistical inverse problem
\eqref{eq:statistical_inverse}, where $X$ is a random function random representing the unknown, and $E$ and $Y$ are the additive noise (statistically independent from $X$) and the observation respectively.

In particular, our prior $\mu_X$ for $X$ is the law of an approximate solution for \eqref{eq:SPDE_inhomogeneous} for some $\alpha > d$, as explained in \Cref{ss:spde-pdos}. Mathematically speaking, $\mu_X$ can be viewed as a Borel measure on the separable Hilbert space $H^{\alpha - d/2 - \epsilon}(\mathbb{T}^d)$ for any $\epsilon > 0$ (see \Cref{ss:whitenoise-sobolev,ss:whitenoise-pdo}). We assume the noise $E$ to be an $m$-variate Gaussian vector with covariance matrix $\Sigma \in \mathbb{R}^{m \times m}$, and the forward operator
\[
    \mathcal{G} \colon H^{\alpha - d/2 - \epsilon}(\mathbb{T}^d) \to \mathbb{R}^{m}
\]
to be continuous. For instance, if $\alpha - \epsilon > d$, the pointwise evaluation operator
\[
    H^{\alpha - d/2 - \epsilon}(\mathbb{T}^d) \owns \xi \mapsto \mathcal{G} (\xi) \mathrel{\mathop:}= \bigl( \xi(x_1), \xi(x_2), \cdots, \xi(x_m) \bigr) \in \mathbb{R}^m
\]
is continuous for any $\{x_1, x_2, \cdots, x_m\} \subset \mathbb{T}^d$ by the Sobolev embedding theorem (see \Cref{ss:whitenoise-pdo}).

With these assumptions, the functional Bayes formula for the inverse problem \eqref{eq:statistical_inverse} reads as
\begin{equation}\label{eq:functional-posterior}
    \frac{\mathrm{d} \mu_{X | Y = \, \boldsymbol{y}}}{\mathrm{d} \mu_X} (\xi) \propto \rho_{\Sigma}\bigl( \boldsymbol{y} - \mathcal{G} (\xi) \bigr), \quad \xi \in H^{\alpha - d/2 - \epsilon}(\mathbb{T}^d),
\end{equation}
where $\rho_{\Sigma} \colon \mathbb{R}^m \to (0, \infty)$ stands for the density function of $E$. We refer to e.g.~\cite[Theorem 6.31]{Stuart_2010} or \cite[Theorem 3.3]{Lasanen-2012-II} for details (and a more general result in case of the latter reference).

Although the distribution of the observation $Y$ is continuous, and so the conditional probability measure $\mu_{X | Y = \, \boldsymbol{y}}$ is strictly defined only for almost all $\boldsymbol{y} \in \mathbb{R}^m$, the right-hand side of \eqref{eq:functional-posterior} is bounded and continuous in both $\boldsymbol{y}$ and $\xi$. This means that we can view the distribution given by \eqref{eq:functional-posterior} as the canonical posterior for \emph{all} $\boldsymbol{y} \in \mathbb{R}^m$, and we have, by the dominated convergence theorem, that
\[
    \mu_{X | Y = \, \boldsymbol{y}'} \longrightarrow \mu_{X | Y = \, \boldsymbol{y}}
    \quad \text{as} \quad
    \boldsymbol{y}' \to \boldsymbol{y}
\]
in the space of Borel probability measures on $H^{\alpha - d/2 - \epsilon}(\mathbb{T}^d)$, e.g.~with repsect to the total variation metric.

\section{Computational Algorithm for Inhomogeneous Priors} \label{sec:computational-method}

In this section we introduce a computational algorithm that samples a random function from the prior distribution $\mu_X$, promoting spatial inhomogeneities. In the interest of readability, the algorithm is presented in the case $d = 1$, where the discretization of the spatial variable $x$ and the frequency variable $\eta$ can be presented as a 2D-tensor. Additional details on the construction of the algorithm in higher dimensions are presented where necessary.

In principle, a draw from the prior $\mu_X$ corresponds to the solution to the SPDE \eqref{eq:SPDE_inhomogeneous} for a right-hand side realization of the white noise. We may look at this solution, formally, as in the discussion following \Cref{eq:SPDE_solution}.

Let us consider a discretization of the 1D torus
\[
    \bar x \mathrel{\mathop:}= \{x_i=i \Delta_x \, | \, i=0, \cdots, N_{x}-1\},
\]
where $\Delta_x = 1/N_{x}$ represents the spatial step size, and a truncated frequency band
\[
    \bar\eta \mathrel{\mathop:}= \{ \eta_i = i - N_{\eta}|i=0,\cdots, 2N_{\eta}\}
\]
for some $N_{\eta} \in \mathbb{N}$. This choice of truncation of the frequency variable allows us to consider a KL-type expansion of the white noise $\Psi$ in \eqref{eq:SPDE_solution}.

Let $\psi \colon \mathbb{T} \to \mathbb{R}$ be a truncated expansion of a realization of $\Psi$ along the Fourier basis functions corresponding to the frequencies $[-N_\eta , N_\eta] \cap \mathbb{Z}$, i.e.
\begin{equation}\label{eq:whitenoise-truncation}
    \psi(x) \mathrel{\mathop:}= \sum_{0 \leq i \leq 2N_{\eta}} \psi_{\eta_i} e^{ \boldsymbol{i} 2\pi \eta_i \cdot x }, \quad x \in \mathbb{T},
\end{equation}
where $\psi_{\eta} \mathrel{\mathop:}= \widehat{\Psi}(\eta)$ stand for the Fourier coefficients of the white noise $\Psi$ for $\eta \in \mathbb{Z}$. For the canonical real-valued white noise, the Fourier coefficients can be sampled as $\psi_0 \sim \mathcal{N}(0, 1)$ and $\psi_\eta = ( a + \boldsymbol{i} b )/\sqrt{2}$ for $\eta \in \mathbb{Z} \setminus \{0\}$, where $a$, $b \stackrel{\text{i.i.d}}{\sim} \mathcal{N}(0, 1)$, with the following dependency structure:
\[
    \psi_{-\eta} = \overline{\psi_\eta} \quad \text{for all } \eta \in \mathbb{Z}, \qquad \psi_{\eta'} \perp \!\!\! \perp \psi_{\eta} \quad \text{if } \eta' \neq -\eta;
\]
see \Cref{ap:whitenoise}. We let a tensor $S$ to hold these Fourier coefficients, i.e.~$[S]_i \mathrel{\mathop:}= \psi_i$ for $0 \leq i \leq 2N_\eta$.

In order to state the computational agorithm involving truncated parametrix approximations (see \Cref{ss:spde-pdos}), let us recall the finite-difference operator $\Delta_\eta^{\gamma}$, and the differential operator $D_x^{(\gamma)}$ familiar from periodic Taylor expansions \cite[Section 3.4]{Ruzhansky_Turunen}.

\begin{definition}
    \begin{enumerate}[(i)]
        \item For functions $F \colon \mathbb{Z} \to \mathbb{C}$, define $\Delta_\eta F$ as the function $\eta \mapsto F(\eta + 1) - F(\eta)$. For integers $\gamma > 1$, define
            \[
                \Delta_\eta^{\gamma} \mathrel{\mathop:}= \underbrace{\Delta_\eta\cdots \Delta_\eta}_{\gamma \text{ times} }.
            \]
            For $d \in \mathbb{N} \setminus \{0\}$, functions $F \colon \mathbb{Z}^d \to \mathbb{C}$ and $\gamma \in \mathbb{N}^d$, define
            \[
                \Delta_\eta^{\gamma} \mathrel{\mathop:}= \Delta_{\eta_1}^{\gamma_1} \cdots \Delta_{\eta_d}^{\gamma_d},
            \]
            with the understanding that $\Delta_{\eta_k}^{\gamma_k} F$ acts on the function
            \begin{equation}\label{eq:difference-coordinatewise}
                \eta' \mapsto F(\eta_1, \cdots, \eta_{k-1}, \eta', \eta_{k+1}, \cdots, \eta_d).
            \end{equation}

        \item For sufficiently regular functions $f \colon \mathbb{T} \to \mathbb{C}$ and $\gamma \in \mathbb{N} \setminus \{0\}$, define
            \[
                D^{( \gamma )}_x f \mathrel{\mathop:}= \prod_{j=0}^{\gamma-1} \left( \frac{1}{\boldsymbol{i}2\pi} \frac{\partial}{\partial x} - j \right) f.
            \]
            For $d \in \mathbb{N} \setminus \{0\}$ and $\gamma \in \mathbb{N}^d$, define the differential operator
            \[
                D^{( \gamma )}_x = D^{(\gamma_1)}_{x_1} \cdots D^{(\gamma_d)}_{x_d},
            \]
            on $\mathbb{T}^d$, where $D_{x_k}^{(\gamma_k)}$ is understood coordinate-wise as in \eqref{eq:difference-coordinatewise}.

        \item For sufficiently regular functions $p \colon \mathbb{T}^d \times \mathbb{Z}^d \to \mathbb{C}$ and $\gamma \in \mathbb{N}^d$, $D^{(\gamma)}_x p$ and $\Delta^{(\gamma)}_{\eta} p$ are understood as actions on the spatial and the frequency variables, respectively.
    \end{enumerate}
\end{definition}

Recall from the discussion in \Cref{ss:spde-pdos} that we are interested in SPDE's associated with the symbol
\begin{equation} \label{eq:p-symbol}
    p_{\alpha}(x,\eta) = \left(\sigma(x)+|\eta|^2\right)^{\alpha/2},
\end{equation}
where $\sigma \in C^{\infty}(\mathbb{T}^d)$ is bounded from below and above by strictly positive constants, $|\cdot |$ is the standard $\ell_2$-norm on $\mathbb{R}$ and $\alpha > d = 1$. We define $P_{\alpha}^{\gamma}$ as a discrete tensor that evaluates $D^{(\gamma)}_x p_{\alpha}$ on $\bar x$ and $\bar \eta$, with respect to $x$, i.e.
\[
    [P_{\alpha}^{\gamma}]_{i,j} = D^{(\gamma)}_x p_{\alpha} (x_i, \eta_j), \qquad i \in \{ 0,\dots,N_{x}-1 \}, ~j \in \{ 0,\dots, 2N_\eta \}.
\]

In the sequel, we will also apply the finite-difference operator $\Delta_\eta^\gamma$ on 2D-tensors $Q$, where we note that, in ``Pythonian'' terms, we have
\[
    \Delta_\eta Q = Q_{:, 1:} - Q_{:, :-1},
\]
and so the operator $\Delta_\eta$ reduces the dimension of the tensor $Q$ in the second component by one. Thus, in order to apply finite differences of order $\gamma$, the original tensor $Q$ should be sampled with appropriately many overhead nodes.

We now provide an approximate discretization for the solution to the SPDE \eqref{eq:SPDE_inhomogeneous} to draw a sample from $\mu_X$. However, instead of considering \eqref{eq:SPDE_solution} directly, we first approximate it by the truncated parameterix
\begin{equation} \label{eq:approximation}
    \xi = [p_{\alpha}(\cdot,D)]^{-1} \Psi \approx q_{-\alpha}^{(N)}(\cdot , D) \Psi,
\end{equation}
where $\xi\sim \mu_X$ and the approximation is understood in the sense of \eqref{eq:approx-error}. By the definition of pseudo-differential operators \eqref{eq:pdo-definition}, we thus have the approximate truncated expansion
\begin{equation}\label{eq:pde-sol-truncation}
    \xi(x) \approx \sum_{0 \leq i \leq 2N_{\eta}} q_{-\alpha}^{(N)}(x,\eta_i ) \psi_{\eta_i} e^{\boldsymbol{i} 2\pi \eta_i \cdot x}, \quad x \in \mathbb{T},
\end{equation}
where the $\psi_{\eta_i}$ are the Fourier coefficients of the white noise \eqref{eq:whitenoise-truncation}.

Let us then define the tensors $q_{-\alpha}^{(N)}(x_i,\eta_j)$, containing the discretized terms in the expansion  \eqref{eq:pde-sol-truncation}. By \cite[Theorem 4.9.13]{Ruzhansky_Turunen}, we have 
\[
    q_{-\alpha}^{(N)}(x_i,\eta_j) = \sum_{k=0}^{N} [Q^k_{-\alpha}]_{i,j},
\]
where summand tensors can be computed following the recursive relation
\begin{align*}
    Q^{0}_{-\alpha} & = [P^{0}_{\alpha}]^{\circ-1}, \\
    Q^{N}_{-\alpha} & = - Q^{0}_{-\alpha} \circ\cdot \left( \sum_{ \substack{ 0 \leq k < N \\ |\gamma| = N - k } } \frac{1}{\gamma!} \Delta_{\eta}^{\gamma}  Q^{k}_{-\alpha} \circ \cdot P^{\gamma}_{\alpha} \right).
\end{align*}
Here $\circ \cdot$ and $\circ-1$ refers to element-wise multiplication and inversion, respectively.

Note that in the special case of spatial homogeneity, i.e.~when $\sigma$ is constant, $P^{\gamma}_{\alpha} \equiv 0$ for $\gamma > 0$, and so $Q^N_{-\alpha} = [P^0_{\alpha}]^{\circ-1}$ for all $N$. In other words, the pseudo-differential operator $p_{\alpha}(\cdot, D)$ is invertible simply by inverting the symbol.

Denote by $\mathcal F^{-1}$ the discrete inverse Fourier transform, and recall that $S$ is the tensor that holds the white noise discretization, i.e.~$[S]_i = \psi_{\eta_i}$ for $0 \leq i \leq 2N_{\eta}$. Then, by \eqref{eq:whitenoise-truncation} and \eqref{eq:pde-sol-truncation}, a discretization of $\xi\sim \mu_X$ can be obtained by
\begin{equation}
    \xi(x_i) \approx \mathcal F^{-1}( [Q^{N}_{-\alpha}]_{i,:} \circ \cdot S ), \qquad i=1,\dots, N_x.
\end{equation}
Here $[Q^{(N)}_{-\alpha}]_{i,:}$ refers to extracting the $i$th slice from the first index.

We summarize the sampling process of $\xi\sim \mu_X$ in the algorithm below.

\begin{algorithm}
    \caption{Sampling $\xi\sim \mu_X$ from an inhomogeneous prior}\label{alg:psido-alg}
    \begin{algorithmic}
        \Require Realization of white noise $S$, symbol $p_\alpha(x,\eta)$ as in \eqref{eq:p-symbol} and truncation term $N$
        \State Compute tensors $P_\alpha^{\gamma}$, for $\gamma=0,\dots,N-1$.
        \State Assemble differentiation tensor $\Delta_\eta$
        \State Compute the base parametrix term $Q^0_{-\alpha}$

        \For{$k=1,\dots,N$}
        \For{$\gamma = 1,\dots,k$}
        \State $\Delta^{\gamma}_{\eta} Q^{k-\gamma} = \Delta_{\eta}[ \Delta^{\gamma-1}_{\eta} Q^{k-\gamma}_{-\alpha} ]$
        \State $T_\gamma := \Delta^{\gamma}_{\eta} Q^{k-\gamma} \circ \cdot P^{\gamma}_{\alpha}$ 
        \EndFor
        \State $Q^k_{-\alpha} = Q^{k-1}_{-\alpha} - Q^{0}_{-\alpha} \circ\cdot \left( \sum_{0 \leq \gamma< k} \frac{1}{\gamma!} T_\gamma \right)$
        \EndFor
        \For{$i=1,\dots,N_x$}
        \State $\xi(x_i) = \mathcal F^{-1}([Q^N_{-\alpha}]_{i,:} \circ \cdot S)$  \Comment{$\mathcal F^{-1}$ is the discrete inverse Fourier transform}
        \EndFor

        \Ensure discrete realization $\xi \sim \mu_X$
    \end{algorithmic}
\end{algorithm}

\section{Numerical Demonstrations} \label{eq:numerical-results}

In this section, we examine the computational performance of the inhomogeneous prior introduced in this paper. We begin with a one-dimensional example, where we numerically investigate the truncation of the pseudo-differential operator $q^{(N)}_{-\alpha}$, and study its realizations. The resulting approximate solutions of the SPDE \eqref{eq:SPDE_inhomogeneous} are compared against a finite-difference discretization of the underlying SPDE.

We then employ this prior in a Bayesian de-noising setting: specifically, we formulate the de-noising problem using an inhomogeneous prior and analyze the resulting posterior distribution. We conclude the numerical experiments with an application to an X-ray computed tomography problem.

\subsection{Prior Discretization} \label{sec:prior}
Here we investigate an inhomogeneous prior of the form \eqref{eq:SPDE_inhomogeneous} with principal symbol \eqref{eq:p-symbol} with components
\begin{equation} \label{eq:psido-symbol}
    \sigma(x) = 0.05 + 2\exp\left( - \frac{(x-0.5)^2}{0.5} \right), \qquad \alpha = 2.
\end{equation}
We first visualize the expansion coefficients $q^{(N)}_{-\alpha}(x,\eta_i)$ in \eqref{eq:pde-sol-truncation}. We begin by discretizing the spatial domain $[0,1]$ into 65 grid points and selecting integer frequencies $-32 \leq \eta \leq 31$. We then follow the first loop in \Cref{alg:psido-alg} to compute the first four terms of the parametrix expansion.

\begin{figure}
    \centering
    \includegraphics[width=\linewidth]{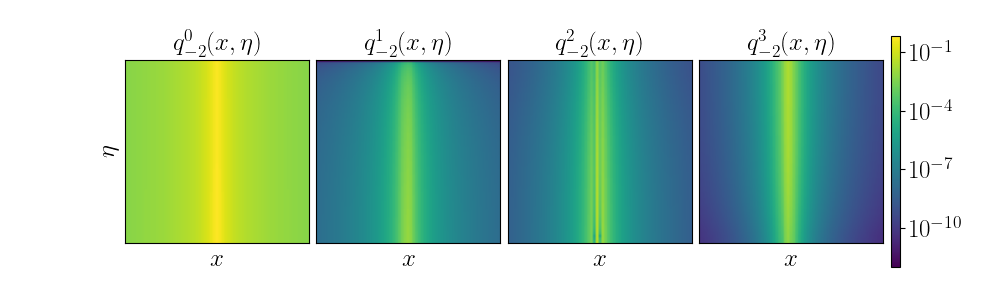}
    \caption{Coefficients of the SPDE \eqref{eq:SPDE_inhomogeneous} with principal symbol \eqref{eq:psido-symbol} }
    \label{fig:q-terms}
\end{figure}

\begin{figure}
    \centering
    \includegraphics[width=0.5\linewidth]{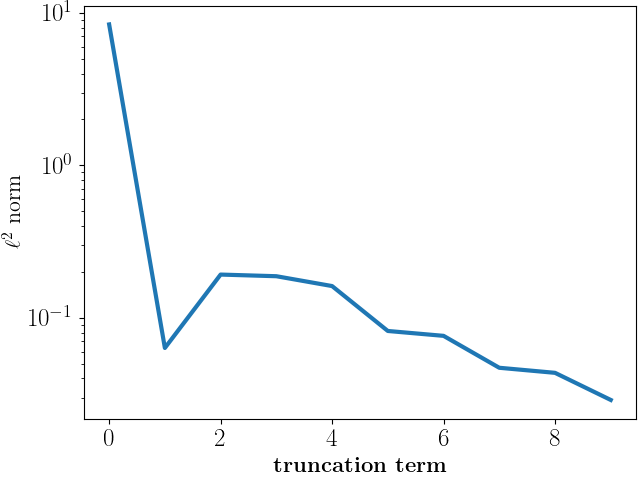}
    \caption{$\ell^2$-norm of the expansion coefficients of the SPDE \eqref{eq:SPDE_inhomogeneous} with principal symbol \eqref{eq:psido-symbol}, up to integration constant.}
    \label{fig:q-terms-norms}
\end{figure}

The expansion coefficients $q^{i}_{-2}(x,\eta)$ are plotted in \Cref{fig:q-terms}. Note that the colormap is displayed on a logarithmic scale, illuminating how rapidly the expansion coefficients diminish. Furthermore, the $0$th term contains most of the mass. This is also demonstrated in \Cref{fig:q-terms-norms}, where we show the norms, up to an integration constant, for the first 10 terms.

These results confirm our visual observation of the rapid decay of the expansion terms -- a fact that is also quantified in \cite[Theorem 4.9.13]{Ruzhansky_Turunen}, in the sense that $q^N_{-\alpha}(\cdot, \cdot) \in S^{-\alpha-N}(\mathbb{T}^d \times \mathbb{Z}^d)$. We note however that the numerical method becomes unstable for large truncation indices. Computational methods for pseudo-differential operators remain an active research area, and developing stable schemes capable of handling large truncations is left for future work.

In the next experiment we investigate the solution of the inhomogeneous SPDE. We create a noise right-hand side for \eqref{eq:SPDE_inhomogeneous} following the truncated white noise \eqref{eq:whitenoise-truncation}.

For comparison, we also consider a finite-difference implementation of \eqref{eq:SPDE_inhomogeneous} by considering
a discretization of $\sigma(x)$ and $\Delta$ with a centered FD method on the real line. The same truncated white noise realization is used as the source term for the two methods.

\begin{figure}
    \centering
    \includegraphics[width=\linewidth]{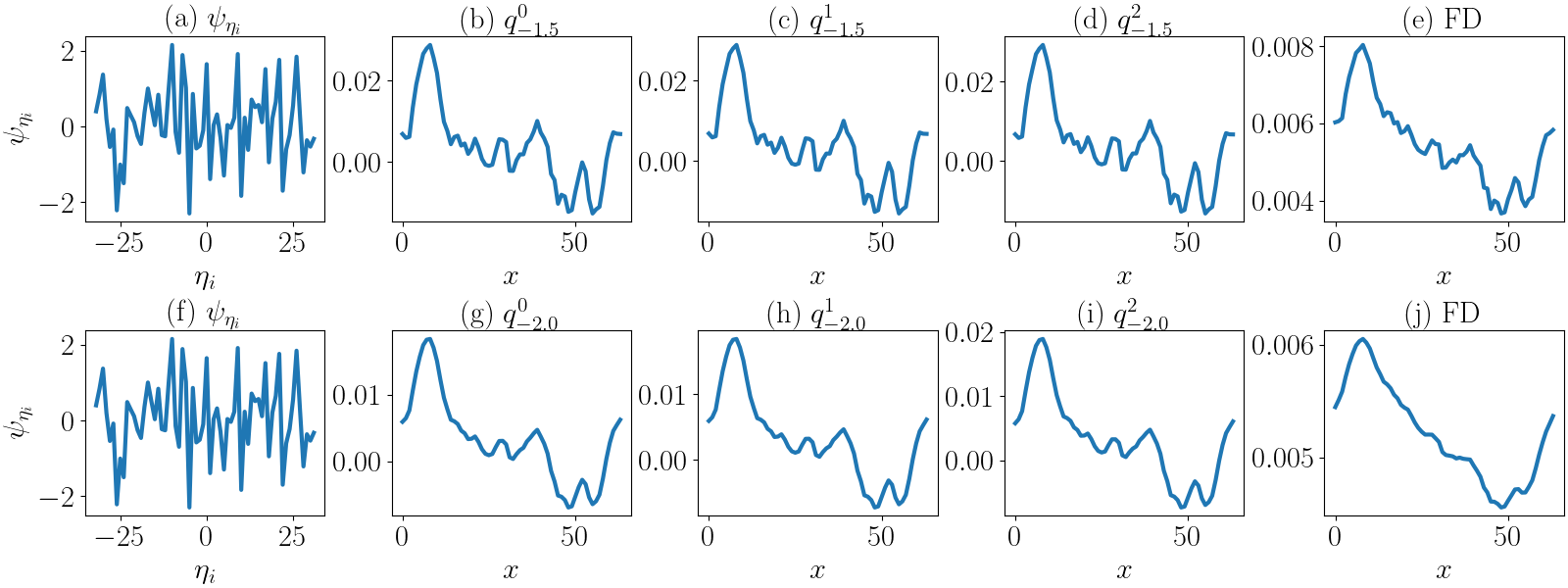}
    \caption{Comparison between the truncated parametrix solution of the SPDE \eqref{eq:SPDE_inhomogeneous} and the finite-difference method.}
    \label{fig:compare-FD}
\end{figure}

\Cref{fig:compare-FD} shows the solution to the inhomogeneous SPDE \eqref{eq:SPDE_inhomogeneous}, following \Cref{alg:psido-alg}, for $\alpha = 1.5$ (first row) and $\alpha = 2$ (second row), using the same white noise realization. Panels (b)--(d) and (g)--(i) correspond to $q^0_{-\alpha}$, $q^1_{-\alpha}$ and $q^2_{-\alpha}$, respectively. The differences among these truncation levels are visually negligible. For comparison, we also provide the FD solutions in panels (e) and (j) for $\alpha = 1.5$ and $\alpha = 2$, respectively.

Recall that $q^{(N)}_{-\alpha}(\cdot, D)\Psi$ and the exact solution to \eqref{eq:SPDE_inhomogeneous} differ by a smoothing term as indicated in \eqref{eq:approx-error}. Thus, although a perfect match in shape or amplitude between the two methods is not expected, the solutions in \Cref{fig:compare-FD} demonstrate a striking similarity. In particular, the roughness characteristics are well captured by the truncated parametrix approximation.

We note that the solution becomes smoother as $\alpha$ increases, which is consistent with the fact that $\alpha$ determines the regularity of the field. However, the discrepancy between the truncated parametrix and the FD approximations grows with increasing $\alpha$, since the smoothing remainder term in \eqref{eq:approx-error} becomes more dominant for larger values of $\alpha$.

\subsection{1D De-noising} \label{sec:de-noising}

Here we employ the inhomogeneous prior introduced earlier within a Bayesian formulation of the de-noising inverse problem. We choose this simple example to illustrate how such priors can be used in practice, and how the results can be interpreted.

The de-noising problem corresponds to \eqref{eq:inverse-summary}, where the forward operator $\mathcal G$ is a uniform spatial sampling of $X$, i.e., an identity operator on a grid. To generate an observation vector $\boldsymbol{y}$, we first draw a realization from the prior, $\xi \sim \mu_X$, and then apply the forward model in \eqref{eq:inverse-summary} to obtain the corresponding realization of $\boldsymbol{y}$. This method of data generation is often referred to as an \emph{inverse crime} \cite{kaipio2005statistical} and is typically avoided when promoting new methods. However, we adopt it here purely for demonstration and reproducibility. We later examine the performance of the approach on an out-of-prior example in two dimensions, where the influence of inhomogeneity is more pronounced.

We first discretize the problem. The discretization of the prior follows the discussion in \Cref{sec:prior} by uniformly discretizing the interval $[0,1]$ into 65 points and introducing integer frequencies $-32 \leq \eta \leq 31$. Furthermore, we truncate the expansion of prior's pseudo-differential operator with $N=5$ terms and follow \Cref{alg:psido-alg} to compute a sample from the prior distribution, for some noise vector $S$.

By denoting $\mathcal P: S\mapsto \xi$ we can introduce the discrete log-posterior as
\begin{equation} \label{eq:NLP}
    \log p(\boldsymbol{\xi|y }) = - \frac{1}{2\sigma^2_{\text{noise} }}\| \boldsymbol{y} - \mathcal P(S) \|^2_2 - \frac{1}{2}\| S \|_2^2 + c,
\end{equation}
where the noise standard deviation $\sigma_\text{noise}$ is chosen following the relation $\sigma_\text{noise} = 0.01\times \| \boldsymbol{y} \|_2$ for a noise free signal $\boldsymbol{y}$, and $c \in \mathbb{R}$ is a constant.

We can explore this posterior using point estimates, such as the maximum a posteriori (MAP) estimator -- i.e.~the $S$ that maximizes \eqref{eq:NLP} -- or the posterior mean. Moreover, we can estimate moments of the posterior distribution to quantify uncertainty. A standard approach for approximating these moments is Markov chain Monte Carlo (MCMC) sampling \cite{kaipio2005statistical}, which generates a point cloud approximating the posterior, with samples concentrated in regions of high posterior probability. Classical examples include the Metropolis–Hastings algorithm; see \cite{kaipio2005statistical,hammersley2013monte,mcbook} for a range of such methods. There also exist MCMC schemes specifically designed for infinite-dimensional Bayesian inverse problems, e.g.~\cite{cotter2013mcmc,boehl2022ensemble}.

In this work, we employ the no-U-turn sampler (NUTS), an adaptive variant of Hamiltonian Monte Carlo (HMC) \cite{hoffman2014no,marwala2023hamiltonian}. HMC augments the parameter space with auxiliary momentum variables and interprets the resulting joint distribution through Hamiltonian dynamics. Proposed samples are generated by approximately simulating trajectories of constant Hamiltonian energy, for which efficient numerical integrators (e.g.~leapfrog schemes) are available. NUTS automatically adapts the trajectory length by detecting when the simulated path begins to double back on itself (a ``U-turn''), thereby removing the need to manually tune the number of leapfrog steps and improving sampling efficiency. While NUTS is not specifically designed for infinite-dimensional posterior distributions, our numerical experiments suggest that it performs robustly in the presence of inhomogeneous priors. Nevertheless, further research is needed to develop sampling strategies that are optimal for the infinite-dimensional setting. In this test problem, we choose 200 samples to warm-up (tune) NUTS method and then draw 2000 samples from the posterior.

We use the following empirical expressions for the posterior mean and the MAP estimate:
\begin{equation}
    S_{\text{mean}} = \frac{1}{N_{\text{sample}}} \sum_{j=1}^{ N_{\text{sample}} } S^{(j)}, \qquad
    S_{\text{MAP}} = \underset{S}{\text{arg max}} \, \log p(\boldsymbol{\xi}(S)|,\boldsymbol{y}),
\end{equation}
and we compute the pointwise $95\%$ highest posterior density interval (HPD) \cite{bernardo1994bayesian} as a measure of uncertainty. Recall that the HPD represents the smallest interval containing $95\%$ of the (pointwise) posterior mass. We denote by $\xi_{\text{mean}}$ and $\xi_{\text{MAP}}$ the push-forwards of $S_{\text{mean}}$ and $S_{\text{MAP}}$ under the mapping $S \mapsto \mathcal{P}(S)$, yielding the corresponding function estimates.

\begin{figure}
    \centering
    \includegraphics[width=\linewidth]{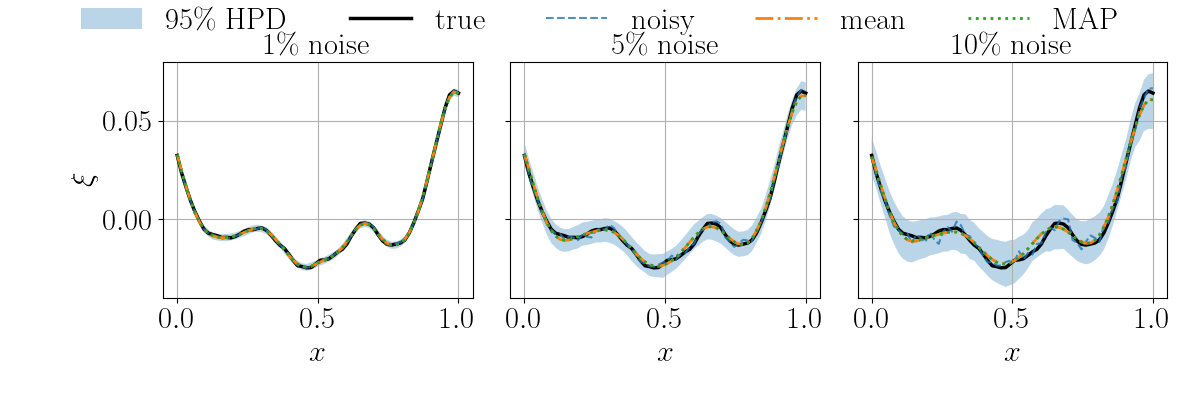}
    \caption{MAP estimation and HPD interval in the denoising problem with various noise levels.}
    \label{fig:denoising-function}
\end{figure}

\begin{figure}
    \centering
    \includegraphics[width=\linewidth]{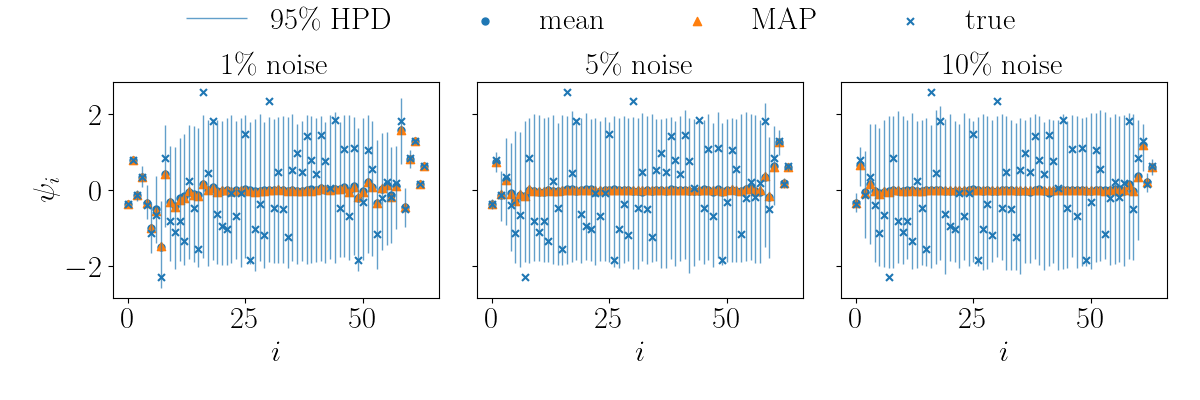}
    \caption{MAP estimation and HPD interval for Fourier coefficients in the denoising problem with various noise levels.}
    \label{fig:denoising-parameters}
\end{figure}

The results for this inference problem are shown in \Cref{fig:denoising-function} for noise levels of $1\%$, $5\%$, and $10\%$. Across all noise levels, we observe that both the posterior mean and the MAP estimate yield continuous functions. For low noise levels, these estimates closely coincide with the true function, whereas for higher noise levels the posterior mean becomes less reliable compared with the MAP estimate. In addition, the uncertainty in the reconstruction increases systematically as the noise level grows.

Furthermore, \Cref{fig:denoising-parameters} presents the estimated coefficients $\psi_i$. We observe that the lower-frequency components are accurately recovered by the method. Additionally, the number of parameters estimated with high certainty decreases as the noise level increases.

To conclude, our results indicate that inhomogeneous priors can be employed in much the same way as their homogeneous counterparts, and can be interpreted analogously to other infinite-dimensional homogeneous priors. Consequently, these methods can be readily integrated into applications where homogeneous priors are commonly used, enabling rapid adoption.

In the next section, we demonstrate the advantages of inhomogeneous priors in settings where the inhomogeneity, i.e.~$\sigma(x)$, is unknown and must be inferred jointly with the function $\xi$. This highlights their applicability beyond what is possible with homogeneous priors.

\subsection{2D Inverse Problems}
In this section, we investigate the performance of inhomogeneous priors in two dimensions and apply them to a limited-angle X-ray computed tomography problem. Our numerical experiments in \Cref{sec:prior} indicate that the parametrix expansion of the pseudo-differential operator for the inhomogeneous prior can be truncated aggressively. Accordingly, we retain only the first term of the expansion (i.e.~the $0$th term), which simplifies the numerical implementation while preserving sufficient flexibility in the prior. We further employ these priors within a hierarchical Bayesian framework in which the inhomogeneity is treated as unknown. The section concludes with an application to limited-angle X-ray CT reconstruction with unknown inhomogeneity.

\subsection{2D Hierarchical Inhomogeneous Priors} \label{sec:hierarchical}

In this section we construct a hierarchical prior distribution that can adapt to inhomogeneities. We first let $\sigma(x)$ be a random function following a $\log$-Gaussian distribution. We can then combine this with in-homogenous prior, discussed earlier in this paper, to create a family of distributions $\mu_X^{\sigma}$, parameterized by their inhomogeneities.

We can formally model such a hierarchical problem as
\begin{equation} \label{eq:hierarchical-prior}
    \begin{aligned}
        Z &\sim \mathcal N(0, a_1( a_2 - \Delta)^4 ), \\
        X &\sim \mu_{X|Z=\sigma(x)}, \quad \text{i.e. } c(x)(10^{(a_3 + \sigma(x))} - \Delta) X = \Psi .
    \end{aligned}
\end{equation}
Here, $Z$ is a random function describing the local inhomogeneities of the samples, so that $10^{(a_3 + Z(x))}$ represents the spatially varying inverse correlation length of realizations of $X$. Here, with an abuse of notation, we use the symbol $\sigma$ to refer to the exponent of the correlation length. The covariance of $Z$, given by $(a_2 - c_{1}\Delta)^4$, produces smooth fields with relatively large spatial correlation lengths, while $a_1$ controls the overall variance. In this test problem we choose $a_2 = 6.25$ and $a_2=2.5$ to produce the desired range in spatial correlation length, and, we draw a realization $\sigma(x)$ from $\mathcal N(0, a_1( a_2 - \Delta)^4 )$ using a Fourier spectral method by solving $( a_2 - \Delta)^2 \sigma(x) = \psi(x)$ for a white noise realization $\psi(x)$ defined in \eqref{eq:whitenoise-truncation} with Fourier coefficients collected in tensor $S_1$. In the Fourier domain this corresponds to
\begin{equation}
    \widehat{\sigma}_{\eta_i} = \frac{1}{c_1} ( 6.25 + |\eta_i|^2 )^{-2} \psi_{\eta_i} ,
\end{equation}
where $\eta_i$ denotes the $i$th 2D multi-index of discrete integer frequencies (e.g.~$\eta_i = (i_1,i_2)$) and $\widehat{\sigma}_{\eta_i}$ is the corresponding discrete Fourier coefficient. Let us collect components $\frac{1}{a_1} ( 6.25 + |\eta_i|^2 )^{-2}$ into a tensor $\Sigma_{\text{hyper}}$. The function $\sigma(x)$ is then obtained by applying the inverse discrete Fourier transform, i.e.~$\mathcal F^{-1}(\Sigma_{\text{hyper}} \circ \cdot S_1)$, which gives us discrete values $\sigma(x_i)$, $i=1,\dots,N_x$. We select $c_1$ such that $\text{Var}(Z)=1$, i.e.
\begin{equation} \label{eq:variance-hyper}
    a_1 = \sum_{|\eta_i| \in \mathbb N^2} (a_2 + |\eta_i|^2 )^{-4},
\end{equation}
which in our numerical setup yields $a_1 = 6.4\times 10^6$.

Once a realization $\sigma(x)$ is generated, we can assemble the pseudo-differential operator corresponding to the covariance of $X$. We apply the approximation in \eqref{eq:approximation} and truncate the symbol expansion after the first term, i.e.
\begin{equation}
    \xi = [p_1(\cdot,D)]^{-1} \Psi \approx q_{-1}^{0}(\cdot,D)\Psi.
\end{equation}
The approximation error in this truncation is characterized by \eqref{eq:approx-error}. Under this simplification, the computation reduces to
\begin{equation} \label{eq:paramterix_2D}
    Q^0_{-1} = [P_{1}^0]^{\circ-1}, \qquad [P_1^0]_{i,j} = p_1(x_i,\eta_j) = c_2(x_i)\bigl(10^{2.5+\sigma(x_i)} + |\eta_i|^2\bigr),
\end{equation}
which provides the zeroth-order approximation to the parametrized pseudo-dif\-fe\-ren\-tial operator.

\begin{figure}
    \centering
    \includegraphics[width=\linewidth]{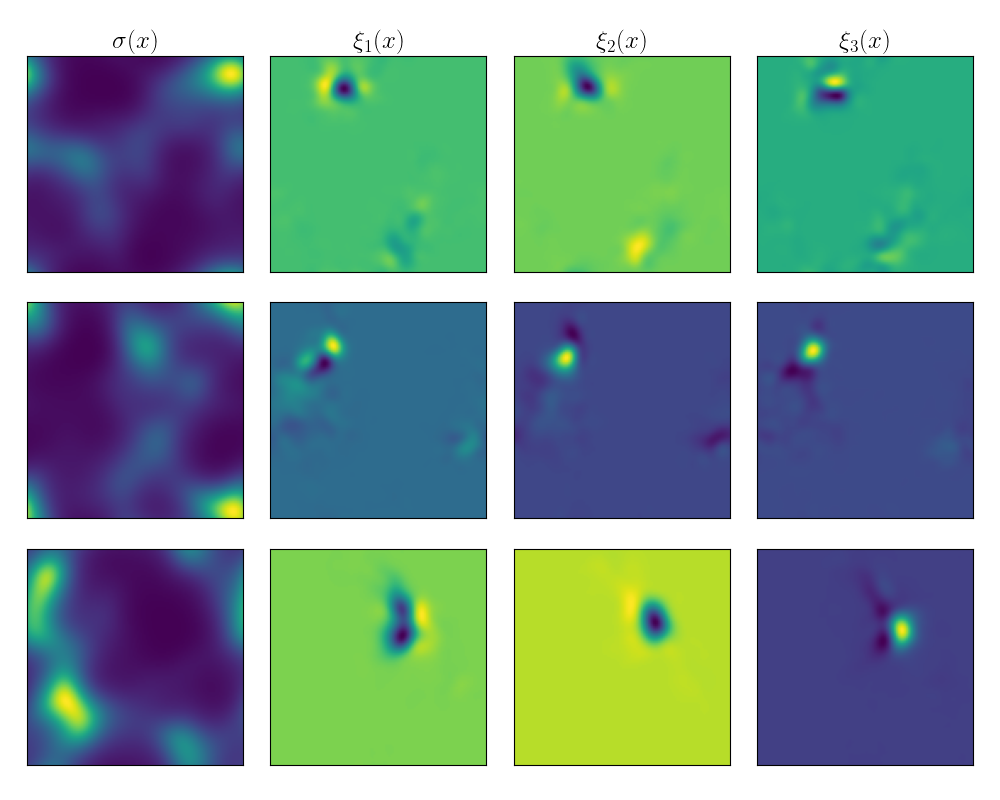}
    \caption{Smooth un-normalized prior samples with inhomogeneity}
    \label{fig:smooth_unnormalized_2D}
\end{figure}

Compared to what is discussed in \Cref{sec:prior}, we introduce the local normalization factor $c(x)$. This is necessary because the spatially varying correlation length $10^{a_3+\sigma(x_i)}$ affects the average amplitude of realizations of $X$. Without normalization, regions with small spatial correlation length would yield realizations with vanishing local variance. This effect is illustrated in \Cref{fig:smooth_unnormalized_2D} and elaborated further below. To prevent this degeneracy, we choose $c(x)$ such that
\begin{equation} \label{eq:normalization_constants}
    c(x_i)^2 = \sum_{\eta_i \in \mathbb N^2 } (10^{2.5+\sigma(x_i)}+|\eta_i|^2)^{-2}.
\end{equation}
The second equation in \eqref{eq:hierarchical-prior} can thus be thought of in terms of solving an SPDE of the form \eqref{eq:SPDE_inhomogeneous} with the colored noise $c(x)^{-1} \Psi$. This can be done exactly as in \Cref{sec:computational-method}, assuming the Fourier transform of this colored noise can be calculated, or approximated to a sufficient extent.

Given a realization of white noise $\psi(x)$ defined in \eqref{eq:whitenoise-truncation}, with its discrete Fourier coefficients collected in a tensor $S_2$, we can compute an approximate sample $\xi$ via
\begin{equation} \label{eq:function_sample}
    \xi(x_i) \approx \mathcal F^{-1}( [Q^0_{-1}]_{i,:} \circ \cdot S_2 ), \qquad i=1,\dots, N_x.
\end{equation}
We summarize the sampling process from this hierarchical prior in \Cref{alg:2D}.

\begin{algorithm}
    \caption{Drawing a sample from the hierarchical inhomogeneous prior \eqref{eq:hierarchical-prior}}\label{alg:2D}
    \begin{algorithmic}
        \Require Realization of white noise $S_1$ and $S_2$ for the hyper-prior and the prior, respectively. Constants $a_2$ and $a_3$.
        \State Compute $a_1$ according to \eqref{eq:variance-hyper}
        \State Compute the tensor $[\Sigma_{\text{hyper}}]_{i} = \frac{1}{a_1} ( a_2 + |\eta_i|^2 )^{-2}$
        \State Compute $\mathcal F^{-1}(\Sigma_{\text{hyper}} \circ \cdot S_1)$ yielding $\sigma(x_i)$, for $i=1,\dots,N_x$
        \State Compute $c(x_i)$ according to \eqref{eq:normalization_constants}
        \State Compute $Q^0_{-1}$ according to \eqref{eq:paramterix_2D}
        \For{$i=1,\dots,N_x$}
        \State find $\xi(x_i)$ from \eqref{eq:function_sample}
        \EndFor
        \Ensure discrete realization $\xi \sim \mu_{X,Z}$
    \end{algorithmic}
\end{algorithm}

\begin{figure}
    \centering
    \includegraphics[width=\linewidth]{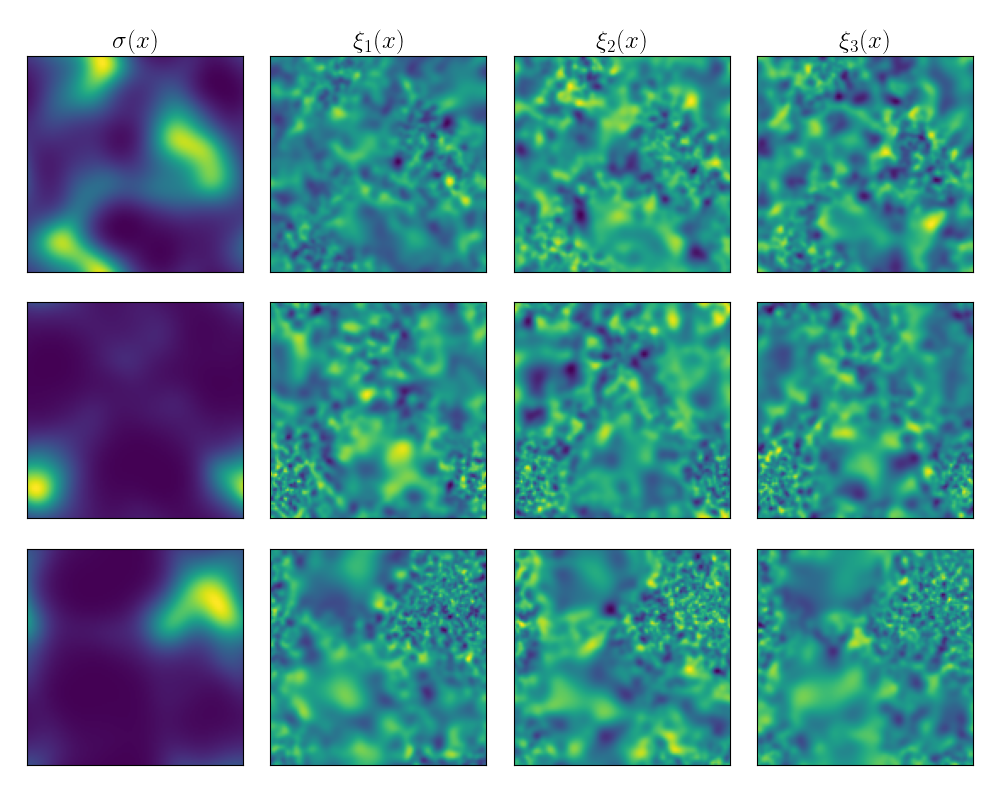}
    \caption{Smooth prior samples with inhomogeneity}
    \label{fig:smooth_inhomogeneoud_2D}
\end{figure}

We provide representative samples from the inhomogeneous hierarchical pri\-or in \Cref{fig:smooth_inhomogeneoud_2D}. The first column displays realizations of the spatially varying correlation length $10^{a_3 + \sigma(x)}$. Fixing this field, we then draw three samples from $X \sim \mu_{X|\sigma(x)}$ as defined in \eqref{eq:hierarchical-prior}. As expected, regions with larger values of $\sigma(x)$ correspond to finer-scale features in $\xi$, reflecting the intended inhomogeneity. For comparison, we also show samples obtained without applying the normalization factor, i.e.~taking $c(x) \equiv 1$ in \eqref{eq:hierarchical-prior} and in \Cref{fig:smooth_unnormalized_2D}. In this case, $\xi$ exhibits pronounced amplitude decay in regions where $\sigma(x)$ is large, demonstrating the vanishing-variance effect discussed earlier.

\begin{figure}
    \centering
    \includegraphics[width=\linewidth]{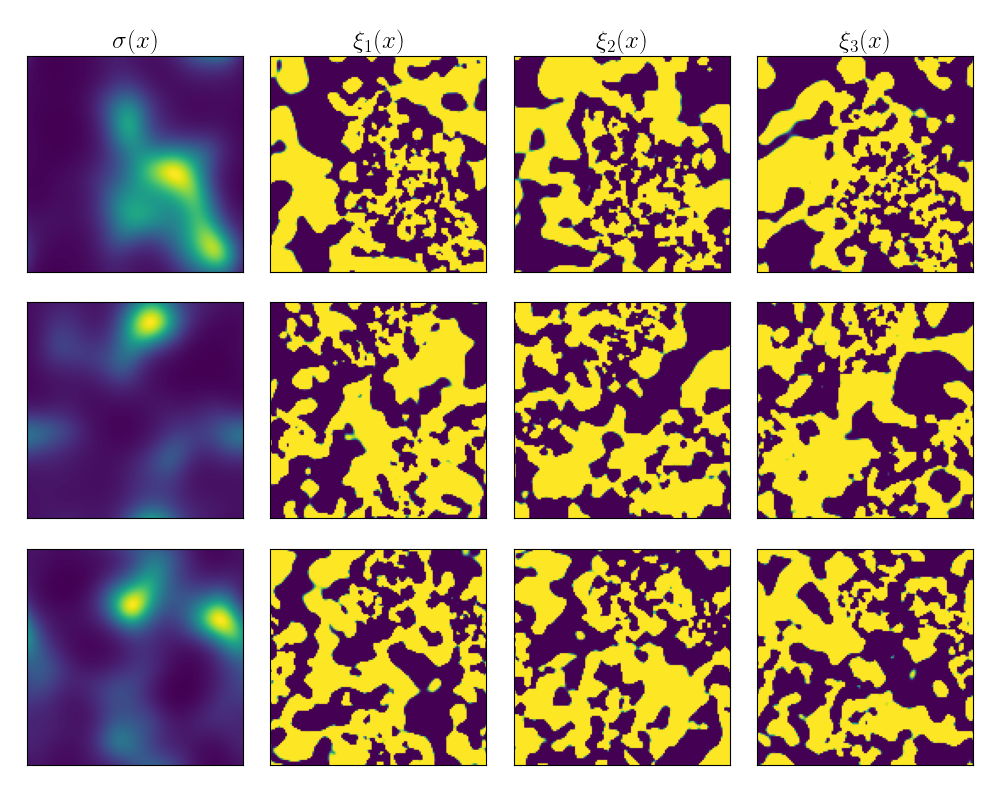}
    \caption{Piecewise constant prior samples with inhomogeneity}
    \label{fig:levelset_2D}
\end{figure}

In many inverse problems, the unknown quantity is often piecewise constant. To accommodate this structure, we define an inhomogeneous level-set prior as a nonlinear push-forward map that transforms $X$ into a piecewise constant random field $\tilde X$. A natural choice would be to use the Heaviside function $H$, setting $\tilde X(x) = c^+ + (c^+ - c^-) H(X)$, where $c^+$ and $c^-$ correspond to the two levels. However, to evaluate the MAP estimation and also utilizing the NUTS sampling algorithm, we require a differentiable transformation that still yields near–piecewise-constant samples. This can be achieved using a sigmoid mapping of the form
\begin{equation}
    \tilde X(x) = \frac{1}{1+ \exp( -kX(x) )},
\end{equation}
where the parameter $k$ controls the sharpness of the transitions. Samples drawn from this inhomogeneous level-set prior are shown in \Cref{fig:levelset_2D}. Similar to the previous case regions with large inverse $\sigma$ correspond to finer inclusions. This prior will be employed in the next section to address a limited-angle X-ray computed tomography problem.

\subsection{Limited Angle X-ray Computed Tomography}

Computed tomography (CT) is a well-established method for reconstructing an image of the cross-section of an object using projection data that represents the intensity loss or attenuation of a beam of X-rays as they pass through the object. According to the Lambert-Beer law \cite{hansen2021computed}, the attenuation of X-rays can be mathematically modeled as the line integrals or projection.

A simple mathematical model of X-ray CT is represented by the \emph{Radon transformation}, represented by a linear operator $\mathcal R$, of an attenuation field, which transforms a 2D attenuation field to its 1D projections, referred to as a \emph{sinogram}. The case where projections are performed along parallel lines is referred to as a parallel beam geometry for X-ray CT. In this case, Radon transform can be formulated in terms of the integral equation
\begin{equation} \label{eq:Radon}
    \mathcal R [\xi] (\theta, s) = \int_{\ell(\theta,s)} \xi(x,y) \ d \ell, \qquad \theta \in [0,\alpha), ~ s\in[-1,1],
\end{equation}
where the integral is an integration along the line $\ell(\theta,s)$ parameterizes by an angle $\theta$ and an intercept $s$. $\theta\in[0,2\pi)$. We assume that $\xi \equiv 0$ outside the unit disk.

In many applications, a complete revolution of $\theta$ within the interval $[0,2\pi)$ is not possible due to physical restrictions. Cases with $\alpha \ll \pi$ are referred to as \emph{limited angle} X-ray CT and are often ill-posed, i.e.~many configurations of unknown $\xi$ may match data.

In this section we consider a limited angle X-ray CT with $\alpha = \pi/4$ while sometimes presenting comparison with a full-angle X-ray CT with $\alpha = \pi$. We discretize $s\in[-1,1]$ into 64 or 128 detector cells and discretize \eqref{eq:Radon} with a Gauss quadrature rule with $64$ or $128$ points, with respect to the discretization level of the detector cells. Furthermore, we discretize $\theta$ to contain 250 equidistant projections in $[0,\pi]$ and extract a subset in the limited angle case. We can formally write the discrete forward map as 
\begin{equation*}
    Y = \mathcal R X + \varepsilon,
\end{equation*}
where $X$ is the unknown cross-sectional image and $\varepsilon \sim \mathcal N(0,\sigma^2_{\text{noise}}I)$ is a Gaussian additive noise with $1\%$ relative noise level, i.e.~if $\boldsymbol{y}_{\text{true}} = \mathcal R \xi_{\text{true}}$, for some true cross-sectional image $\xi_{\text{true}}$, then $\sigma_{\text{noise}} = 0.01\times\| \boldsymbol{y}_{\text{true}} \|_2$.

To set up the Bayesian inverse problem for this X-ray CT problem, we introduce inhomogenuity hyper-parameter $Z$ and equip $(X,Z)$ with the hierarchical in homogeneous prior \eqref{eq:hierarchical-prior}, and define the joint log-posterior as
\begin{equation} \label{eq:cost-hierarchical}
    \log p(X, Z| Y) = - \frac{1}{2\sigma^2_{\text{noise}}} \| \boldsymbol{y} - \mathcal R(\xi(\sigma) )\|_2^2 - \| S_1 \|^2_2 - \| S_2 \|^2_2 + c.
\end{equation}
where $S_1$ and $S_2$ are discretizations of $\Psi_1$ and $\Psi_2$ discussed in \Cref{sec:hierarchical}, $\boldsymbol{y}$ corresponds to X-ray attenuation data, and by $\xi(\sigma)$ we refer to the concatenation of the maps $S_1\mapsto\sigma(x)$ and $(\sigma(x),S_2)\mapsto \xi$, following \Cref{sec:hierarchical}.

\begin{figure}
    \centering
    \includegraphics[width=\linewidth]{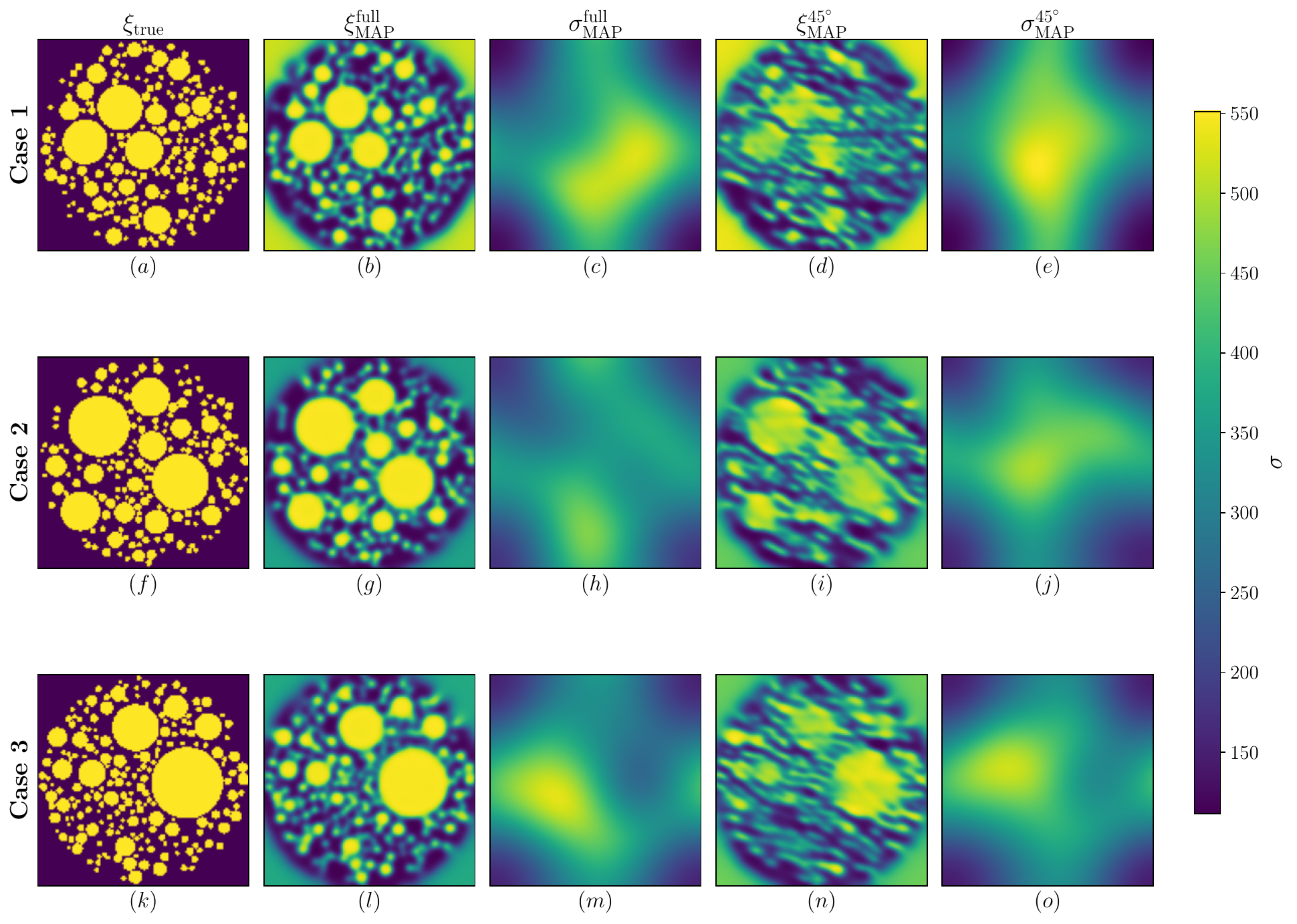}
    \caption{MAP estimation of the attenuation and inhomogeneity of 2 randomly generated and out of prior attenuation field. The first column present the ground truth attenuation fields, the second and third columns correspond to MAP estimation for a full angle X-ray CT, respectively, and the last two columns corresponds to the MAP estimate of the attenuation and inhomogeneity for a limited-angle X-ray CT.}
    \label{fig:MAP-multiple-cases}
\end{figure}

\begin{figure}
    \centering
    \includegraphics[width=\linewidth]{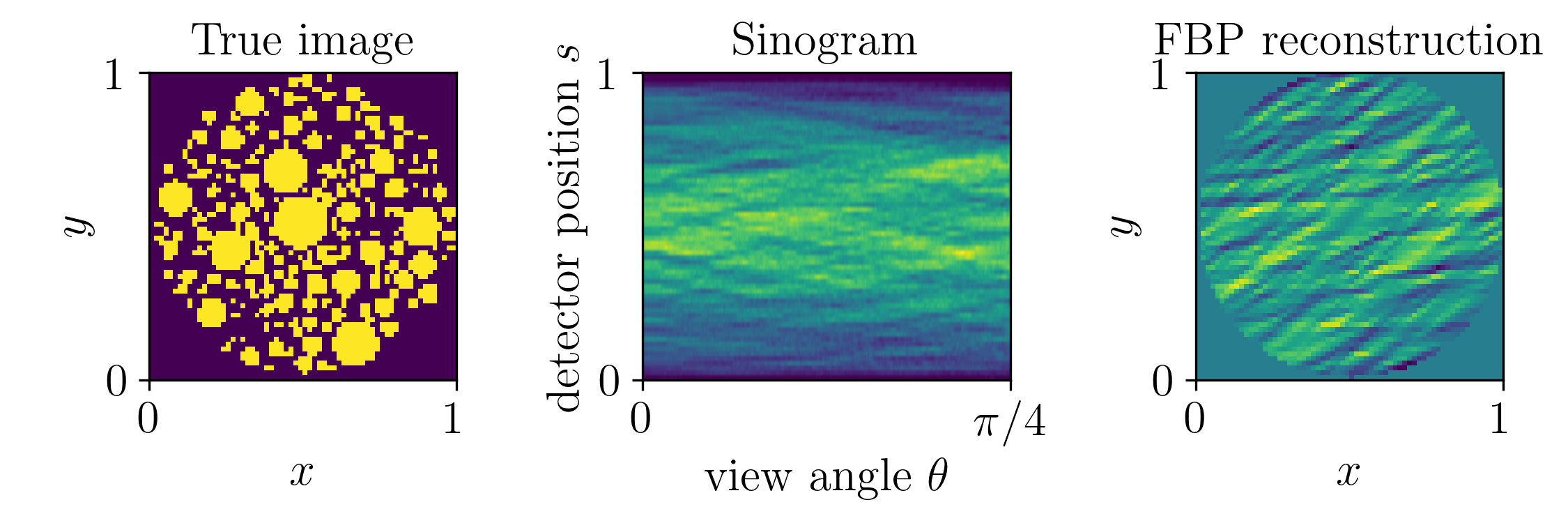}
    \caption{True X-ray image, together with its limited angle attenuation data and an FBP reconstruction.}
    \label{fig:xray-true}
\end{figure}

To create a piece-wise constant and out-of-prior, inhomogeneous sample, we follow a greedy algorithm that places disk inclusions with random diameter and center in the unit disk. In such an algorithm, when non-overlapping of inclusion condition is enforced, the resulting image will naturally contain regions with smaller inclusions and regions with larger inclusions. The ground truth images used in this section are provided in \Cref{fig:MAP-multiple-cases,fig:xray-true}.

We suggest approximately evaluating the MAP point of this hierarchical problem as a point estimate of this limited angle X-ray CT problem. We use the L-BFGS optimization algorithm to maximize the cost function defined in \eqref{eq:cost-hierarchical}. L-BFGS is a quasi-Newton optimization algorithm that approximates the inverse Hessian or the negative log-posterior using only a small number of stored vector pairs, making it highly efficient for large-scale problems \cite{nocedal1980updating}. In practice, we use automatic differentiation and PyTorch \cite{pytorch} to solve the maximization problem.

We estimate the MAP point for 3 randomly generated images with the precedure discussed above. For this experiment we choose 128 discreization points for the detector pixels and the quadrature approximation of line integrals. We repeat the experiment once for a full-angle setup, and once again for a limited angle setup. The results can be found in \Cref{fig:MAP-multiple-cases} with a shared colorbar for the estimated $\sigma$. We report that the Hierarchical prior can correctly identify the regions of small spatial correlation as indicated in the MAP estimation of $\sigma$. Although we notice a degrade in the estimation of $\xi$ and $\sigma$ in a limited angle case, reconstructions are consistent with the full angle reconstruction, as well as, with the ground truth image. In addition, in Case 2, where the phantom exhibits a more uniform distribution of large inclusions, $\sigma$ is estimated with less spatial variation compared to other cases, suggesting that a more homogeneous phantom is detected. We remark that a 45-degree limited angle X-ray CT is a challenging computational problem and our reconstruction significantly improves standard reconstruction approaches (see \Cref{fig:xray-true} for a typical reconstruction). Therefore, we conclude that by identifying the regions of high and low spatial correlation, the hierarchical method has reduced the dependency of the reconstruction to data.

\begin{figure}
    \centering
    \includegraphics[width=\linewidth]{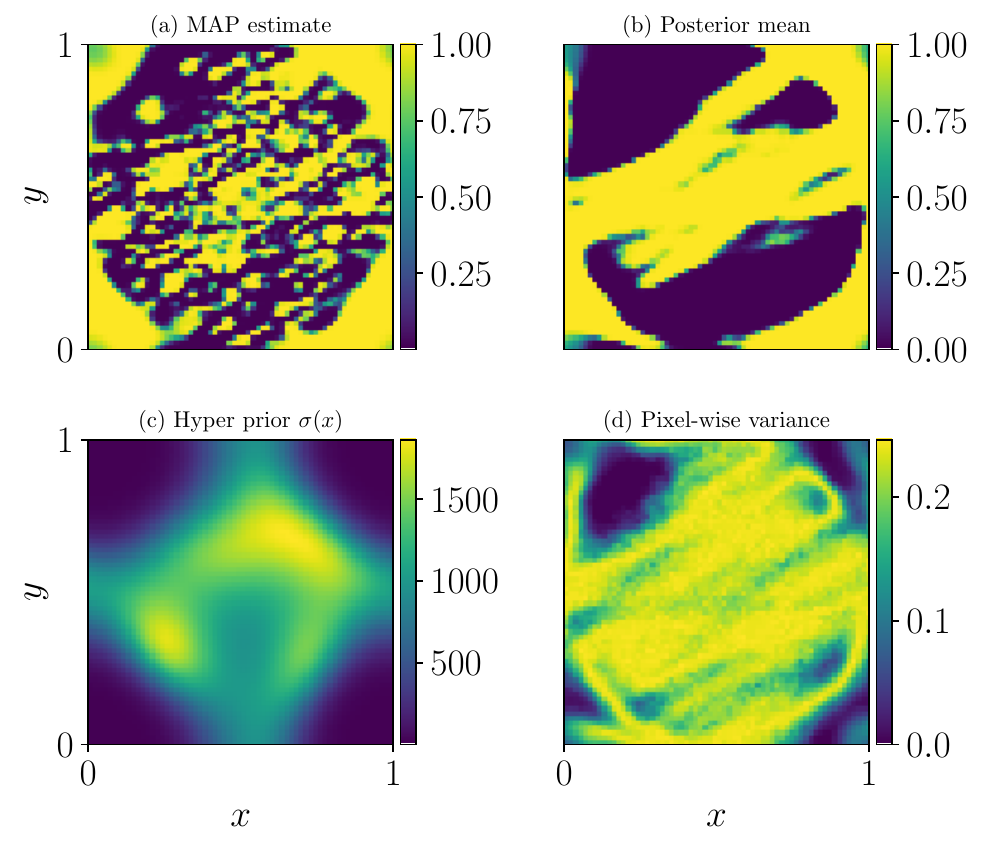}
    \caption{Estimation and uncertainty visualization of the X-ray CT problem with hierarchical inhomogeneous prior. This plot shows (a) the MAP estimation of the image, (b) posterior mean estimation of the image, (c) MAP estimation of the inhomogeneity (d) pixel-wise variance as a measure of uncertainty}
    \label{fig:xray-uq}
\end{figure}

Next, we will fix $\sigma$ to its MAP estimation and attempt to quantify uncertainties in the reconstruction of $\xi$. We create a new random ground truth attenuation field and discretize the problem with 64 detector pixels and quadrature points. MAP estimation of $\sigma(x)$ and $\xi$, for this test problem, are presented in \Cref{fig:xray-uq}. We notice that larger values of $\sigma$ correctly identifies regions with small inclusions. Furthermore, the MAP estimation of also contain smaller inclusions in regions with high value of $\sigma$. We remark that the MAP estimation of $\xi$ prioritizes inclusion, even in a limited angle setting. As a comparison, we provide a filtered back-projection (FBP) \cite{hansen2021computed} reconstruction of this limited angle X-ray CT problem in \Cref{fig:xray-true}. The reconstruction from this method significantly improves the FBP reconstruction. We remark that there are more advanced reconstruction methods prioritizing piece-wise constant images, e.g.~with total variation regularized reconstruction \cite{hansen2021computed}. However, with such data-sparsity (limited angle setup), we do not expect significant improvement compared to the FBP reconstruction. By identifying the regions with small/large correlation length, our method can dramatically reduce the dependency of data for reconstruction.

To quantify uncertainties in the reconstruction, we fix $\sigma(x)$ to its MAP estimation and then perform the UQ procedure for $\xi$ as discussed in \Cref{sec:de-noising}. We use the NUTS sampling method to draw samples from the conditional (conditioned on the MAP estimation of $\sigma(x)$) posterior. We use the NUTS method from Pyro Python package \cite{bingham2018pyro} with 200 warm-up steps and 2000 posterior samples. We report that the effective sample size (computed by the built-in functionality of Pyro which is based on auto-correlation function of the samples) for all components of $S_2$ were at least 1628, indicating high quality of posterior samples for the purpose of the ergodic estimation of the mean and posterior variance.

\begin{figure}
    \centering
    \includegraphics[width=\linewidth]{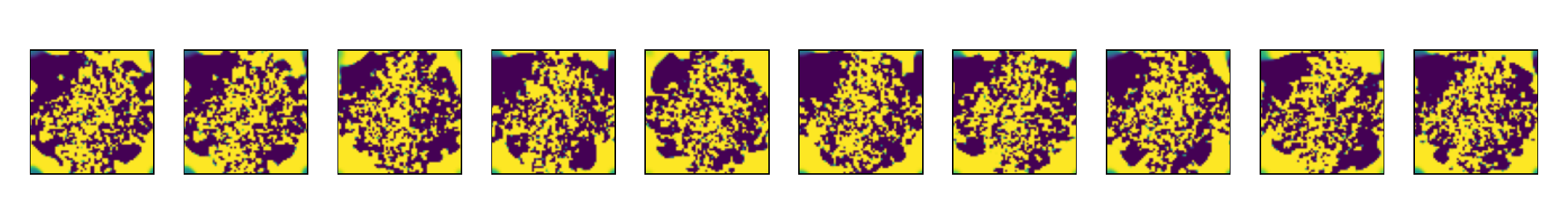}
    \caption{Samples from the posterior distribution of the X-ray CT problem with hierarchical inhomogeneous prior}
    \label{fig:posterior-samples}
\end{figure}

We provide the posterior mean and pointwise variance of $\xi$ in \Cref{fig:xray-true}. We can see a significant difference between the mean and the MAP estimation. We remind the reader that the X-ray CT problem formulated in this section is non-linear, due to the piecewise constant nature of the prior. For a nonlinear and multimodal posterior, the MAP and mean often do not coincide. For comparison, we also provide some posterior samples in \Cref{fig:posterior-samples}. It is clear that posterior samples still preserve the inclusion structure of the prior, similar to the MAP estimation. Therefore, we conclude that the posterior mean may not be an appropriate point estimate for this problem.

We also provide the pixel-wise variance, as an indicator of uncertainty in the estimation. We notice that the orientation of uncertainty is aligned with 45 degree lines confirming expected uncertainty by micro-local analysis \cite{hansen2021computed}. In addition, regions with high uncertainty correspond to some of the artifacts in the MAP estimation.

\section{Conclusion} \label{sec:conclusion}
We introduced a new class of inhomogeneous Gaussian priors for Bayesian inverse problems by defining nonstationary Whittle--Mat\'ern--type fields through convolution with white noise, equivalently characterized as solutions to a pseudo-differential SPDE. The resulting priors fit within the standard infinite-dimensional Bayesian well-posedness framework and admit efficient sampling via a parametrix-based approximation with quantified truncation error. Numerical results in 1D denoising and 2D limited-angle X-ray CT demonstrate improved reconstructions and more informative uncertainty quantification in data-limited regimes.

Natural directions for future work include incorporating locally anisotropic structure by modifying the differentiation/operator symbol in the pseudo-differential formulation, developing dedicated sampling strategies tailored to these priors in large-scale inverse problems, and designing numerically robust schemes that remain stable for higher-order truncations of the parametrix expansion.

\appendix
\crefalias{section}{appendix}
\crefalias{subsection}{appendix}

\section{White Noise Analysis}\label{ap:whitenoise}

In this section we briefly review the existence and some mapping properties of the Gaussian white noise on the $d$-dimensional torus $\mathbb{T}^d \mathrel{\mathop:}= [0,1)^d$. For this purpose, we recall the space $C^{\infty}(\mathbb{T}^d)$ of smooth, complex-valued periodic test functions, equipped with the topology generated by the seminorms
\[
    f \mapsto \| \partial^{\alpha} f \|_{L^{\infty}(\mathbb{T}^d)}, \quad \alpha \in \mathbb{N}_0^d.
\]
We refer to e.g.~\cite[Chapters 3 and 4]{Ruzhansky_Turunen} for a more detailed presentation on periodic analysis and the definitions of several other function spaces appearing here.
We define the Fourier transform $\widehat{f} \colon \mathbb{Z}^d \to \mathbb{C}$ of a function $f \in C^{\infty}(\mathbb{T}^d)$ as
\[
    \mathbb{Z}^d \owns \eta \mapsto \int_{\mathbb{T}^d} f(x) e^{-\boldsymbol{i} 2\pi \eta \cdot x} \, \mathrm{d} x \in \mathbb{C},
\]
where ``$\mathrm{d} x$'' stands for the standard (unnormalized) Lebesgue measure on $\mathbb{T}^d$.

Let $w_0$, $(a_\eta)_{ \eta \in \mathbb{Z}^d \setminus \{0\} }$ and $ (b_\eta)_{ \eta \in \mathbb{Z}^d \setminus \{0\} } $ be independent random variables distributed as $ \mathcal{N}(0, 1) $.
The (real) white noise $\Psi$ on $\mathbb{T}^d$ is defined as
\begin{equation}\label{eq:whitenoise-definition}
    \Psi(x) \mathrel{\mathop:}= w_0 + \sum_{\eta \in \mathbb{Z}^d \setminus \{0\} } \left( a_\eta \cos(x \cdot \eta) + b_\eta \sin(x \cdot \eta) \right),
\end{equation}
where the convergence of the series is understood in the sense of tempered distributions. More precisely, we have
\[
    \left| \int_{\mathbb{T}^d} f(x) \cos(x \cdot \eta) \, \mathrm{d} x \right| \leq C(f, d) (1 + |\eta|)^{-(d+1)}
\]
for all $\eta \in \mathbb{Z}^d$ and test functions $f \in C^{\infty}(\mathbb{T}^d)$ and similarly for $\sin(\cdot)$ in place of $\cos(\cdot)$, so because e.g.~$\mathbb{E} [ \sum_{\eta \in \mathbb{Z}^d \setminus \{0\} } (1+|\eta|)^{-d-1} ( |a_\eta| + |b_\eta| ) ] < \infty$, the right-hand side of \eqref{eq:whitenoise-definition} defines a Borel-measurable random continuous functional of $f \in C^{\infty}(\mathbb{T}^d)$ with full probability.

Relatively simple computations reveal that the Fourier transform of $\Psi$ is given for frequency $\eta = 0$ as $\widehat{\Psi}(0) = w_0$, and for $\eta \in \mathbb{Z} \setminus \{0\} $ by
\begin{equation}\label{eq:whitenoise-fourier}
    \widehat{\Psi}(\eta) = \frac{a_{\eta} + a_{-\eta}}{2} + \boldsymbol{i} \frac{b_{-\eta} - b_{\eta}}{2} =\mathrel{\mathop:} w_{\eta},
\end{equation}
so for $\eta \in \mathbb{Z} \setminus \{0\}$ the Fourier coefficients satisfy $w_{\eta} = \overline{w_{-\eta}}$, and otherwise they are mutually independent circularly contoured complex Gaussian random variables.

Further computations show that $\Psi$ is a real-valued Gaussian field in the sense that
\[
    \left( \langle \Psi, f_1 \rangle, \langle \Psi, f_2 \rangle, \cdots, \langle \Psi, f_n \rangle \right)
\]
is a centered normal random vector for all real-valued test functions $(f_j)_{j \in 1{:}n} \subset C^{\infty}(\mathbb{T}^d)$, with covariance structure
\begin{equation}\label{eq:whitenoise-cov-structure}
    \mathbb{E} \left[ \langle \Psi, f \rangle \langle \Psi, g \rangle \right] = \langle f, g \rangle_{L^2(\mathbb{T}^d)}
\end{equation}
for all real-valued test functions $f$ and $g$. Thus, $\langle \Psi, f \rangle$ and $\langle \Psi, g \rangle$ are independent whenever $fg \equiv 0$. This latter fact means that, although $\Psi$ can not be realized as a random pointwise-defined function, in a distributional sense it has ``independent values at every point'' \cite[Chapter III, Section 4]{gelfand1964generalized4}.

\subsection{White Noise in Sobolev Spaces}\label{ss:whitenoise-sobolev}

Recall that for $s \in \mathbb{R}$, the Sobolev space $H^s(\mathbb{T}^d)$ is defined as the linear space of tempered distributions $u \in \mathcal{D}'(\mathbb{T}^d)$ such that
\[
    \| u \|_{H^s(\mathbb{T}^d)} \mathrel{\mathop:}= \sqrt{ \sum_{\eta \in \mathbb{Z}^d} (1 + |\eta|^2)^s |\widehat{u}(\eta)|^{2} } < \infty.
\]
Note that each summand in the series above is a continuous function of $u \in \mathcal{D}'(\mathbb{T}^d)$, which means that the extended norm
\[
    \| \cdot \|_{H^s(\mathbb{T}^d)} \colon \mathcal{D}'(\mathbb{T}^d) \to [0, \infty]
\]
is a Borel-measurable function on the space of tempered distributions.

Let $\epsilon > 0$, fixed from now. The white noise $\Psi$ is $(-d/2-\epsilon)$-regular in the sense that 
\[
    \mathbb{P} \left( \Psi \in H^{-d/2-\epsilon}(\mathbb{T}^d) \right) = 1.
\]
This result can be found e.g.~in \cite[Theorem 3.4]{veraar2010regularitygaussianwhitenoise} for the complex white noise on $\mathbb{T}^d$, but for our purposes it suffices to note that, by \eqref{eq:whitenoise-fourier},
\[
    \mathbb{E} \left[ \| \Psi \|_{H^{-d/2-\epsilon}(\mathbb{T}^d)}^2 \right] = \sum_{\eta \in \mathbb{Z}^d} (1+|\eta|^2)^{-d/2-\epsilon} \mathbb{E} \left[ |w_\eta|^2 \right] < \infty.
\]

Let us then note that a version of the law of $\Psi$ can be realized as a Borel probability measure on the (separable Hilbert) space $H^{-d/2-\epsilon}(\mathbb{T}^d)$. The argument for this is essentially contained e.g.~in the proof of \cite[Proposition 3.1]{Kusuoka1982TheSP}, but for the sake of completeness we briefly recall the details below. See also \cite[Appendix A]{Helin_Lassas_Oksanen_2014} for a similar result for the Gaussian white noise in local Sobolev spaces over the real line $\mathbb{R}$.

The crux is that for each $f \in H^{-d/2-\varepsilon}(\mathbb{T}^d)$ and $r > 0$, we can write the closed $H^{-d/2-\epsilon}$-ball at $f$ with radius $r$ as
\[
    \bigcap_{N \in {\mathbb{N}}} \left\{ u \in \mathcal{D}'(\mathbb{T}^d) \,:\, \sum_{\eta \in \mathbb{Z}^d, \, |\eta| \leq N} (1 + |\eta|^2)^{-d/2-\epsilon} |\widehat{u - f}|^2 \leq r^2 \right\},
\]
where each set in the intersection is by continuity a closed subset of $\mathcal{D}'(\mathbb{T}^d)$. This means that closed balls in $H^{-d/2-\epsilon}(\mathbb{T}^d)$ are Borel subsets of $\mathcal{D}'(\mathbb{T})$, and this then extends to all Borel subsets of $H^{-d/2-\epsilon}(\mathbb{T}^d)$.

\subsection{White Noise and Pseudo-differential Operators}\label{ss:whitenoise-pdo}

Recall the H\"or\-man\-der symbol class $S^m( \mathbb{T}^d \times \mathbb{Z}^d )$, $m \in \mathbb{R}$, introduced previously in \Cref{ss:whitenoise-pdo}. The action of a pseudo-differential operator $q(\cdot, D)$ with symbol $q(\cdot, \cdot) \in S^m( \mathbb{T}^d \times \mathbb{Z}^d )$ on tempered distributions is defined as follows: if $u \in \mathcal{D}'(\mathbb{T}^d)$, then $q(\cdot, D) u \in \mathcal{D}'(\mathbb{T}^d)$ is defined by
\begin{equation}\label{eq:pdo-transpose}
    \langle q(\cdot, D)u, f \rangle = \langle u, q^\intercal(\cdot, D)f \rangle, \quad f \in C^{\infty}(\mathbb{T}^d),
\end{equation}
where $q^\intercal(\cdot, D)$ is the transpose of the operator $q(\cdot, D)$ (with respect to the dual pairing between $\mathcal{D}'(\mathbb{T}^d)_{|C^{\infty}(\mathbb{T}^d)}$ and $C^{\infty}(\mathbb{T}^d)$) with symbol $q^\intercal(\cdot, \cdot) \in S^m( \mathbb{T}^d, \mathbb{Z}^d )$; see e.g.~\cite[Theorem 4.7.3]{Ruzhansky_Turunen}.
Note that if the symbol $q(\cdot, \cdot)$ satisfies
\begin{equation}\label{eq:symbol-real}
    q(x, \eta) = \overline{q(x, -\eta)} \qquad \forall (x, \eta) \in \mathbb{T}^d \times \mathbb{Z}^d
\end{equation}
(as is the case in our applications), then $q(\cdot, D)$ and $q^\intercal(\cdot, D)$ map real-valued test functions to real-valued test functions.

Recall then from the previous subsection that the white noise $\Psi$ in our setting is a random element in the Sobolev space $H^{-d/2-\epsilon}(\mathbb{T}^d)$, where $\epsilon > 0 $ is arbitrary and fixed.

Let $q(\cdot, D)$ be a pseudo-differential operator in the class $S^{-\alpha}(\mathbb{T}^d \times \mathbb{Z}^d)$ for some $\alpha > 0$, satisfying \eqref{eq:symbol-real}. We have by e.g.~\cite[corollary 4.8.3]{Ruzhansky_Turunen} that $q(\cdot, D) \Psi$ is a random element in the Sobolev space $H^{\alpha-d/2-\epsilon}(\mathbb{T}^d)$. By the discussion surrounding \eqref{eq:whitenoise-cov-structure} and \eqref{eq:pdo-transpose}, this random distribution is a real-valued Gaussian field with the covariance structure
\[
    \mathbb{E} \left[ \langle q(\cdot, D) \Psi, f \rangle \langle q(\cdot, D) \Psi, g \rangle \right] = \left\langle q^\intercal(\cdot, D)f, q^\intercal(\cdot, D)g \right\rangle_{L^2(\mathbb{T}^d)}
\]
for all test functions $f$ and $g$, and so $q(\cdot, D) \Psi$ has the covariance operator
\[
    \mathcal{C} = q(\cdot, D) \circ q^\intercal(\cdot, D)
\]
in the sense of \eqref{eq:gaussian-hilbertspace}.

In particular, if $\alpha > d$, we can take $\epsilon > 0$ so that $\alpha > d + \epsilon$. Then $q(\cdot, D) \Psi$ is a random field with Sobolev regularity of order
\[
    \alpha - d/2 - \epsilon > d/2,
\]
so by the Sobolev embedding theorem, the random field has sample functions in the H\"older space $C^{\alpha - d-\varepsilon}(\mathbb{T}^d)$ (or the Zygmund space of order $\alpha - d - \varepsilon$ in case $\alpha - d - \varepsilon$ is an integer) \cite[Theorems 3.5.4 and 3.5.5]{Schmeisser_Triebel}.

Although the eigenvalues of the covariance operator $\mathcal{C}$ above are generally not easily identified, in this case the decay of $q(\cdot, \cdot) \in S^{-\alpha}( \mathbb{T}^d \times \mathbb{Z}^d )$ with $\alpha > d$ readily implies the following Karhunen-Lo\`eve-like expansion for the random function $q(\cdot, D)\Psi$:
\[
    q(\cdot, D)\Psi(x) = \sum_{ \eta \in \mathbb{Z}^d } q(x, \eta) w_{\eta} e^{\boldsymbol{i} 2\pi \eta \cdot x},
\]
with absolute pointwise convergence with respect to $x \in \mathbb{T}^d$ with full probability.

Convergence rates for the above expansion can be estimated e.g.~as follows. For a symbol $q(\cdot, \cdot)$ satisfying \eqref{eq:symbol-real}, the absolute pointwise convergence in conjunction with Fubini's theorem yields
\begin{align*}
    & \sup_{x \in \mathbb{T}^d} \mathbb{E} \left[ \left| q(\cdot, D)\Psi(x) - \sum_{ |\eta| \leq M } q(x, \eta) w_{\eta} e^{\boldsymbol{i} 2\pi \eta \cdot x} \right|^2 \right]
    \\
    & \qquad = \sup_{x \in \mathbb{T}^d} \qquad \sum_{|\eta| > M} |q(x, \eta)|^2
    \\ 
    & \qquad \leq \sup_{x \in \mathbb{T}^d , \, \eta \in \mathbb{Z}^d } \left( |\eta|^\alpha \left| q(x, \eta) \right| \right)^2 \sum_{|\eta| > M} |\eta|^{-2\alpha}
    \\
    & \qquad \leq \sup_{x \in \mathbb{T}^d , \, \eta \in \mathbb{Z}^d } \left( |\eta|^\alpha \left| q(x, \eta) \right| \right)^2 \sum_{k = M}^\infty A_d(k) k^{-2\alpha},
\end{align*}
where $A_d(k)$ stands for the number of lattice points $\eta \in \mathbb{Z}^d$ with $k < |\eta| \leq k+1$. This number grows like $k^{d-1}$ (and the growth can be quantified more explicitly e.g.~for $d \in \{1, 2\}$), so the the latter series tends to zero like $M^{d - 2\alpha}$ as $M \to \infty$. Let us also note that in our 2D experiments in \Cref{eq:numerical-results}, where $q(\cdot, \cdot) \mathrel{\mathop:}= q_{-\alpha}(\cdot, \cdot)$ is simply given as $p_{\alpha}(\cdot, \cdot)^{-1}$, the latter supremum is $\leq 1$.


\bibliographystyle{siamplain}
\bibliography{references}
\end{document}

%% file: ex_shared.tex
\usepackage{amsfonts}
\usepackage{graphicx}
\usepackage{epstopdf}
\usepackage{amssymb}
\usepackage{enumerate}
\usepackage{algorithm}
\usepackage{algpseudocode}

\ifpdf
  \DeclareGraphicsExtensions{.eps,.pdf,.png,.jpg}
\else
  \DeclareGraphicsExtensions{.eps}
\fi

\newsiamremark{remark}{Remark}
\newsiamremark{hypothesis}{Hypothesis}
\crefname{hypothesis}{Hypothesis}{Hypotheses}
\newsiamthm{claim}{Claim}
\newsiamremark{fact}{Fact}
\crefname{fact}{Fact}{Facts}

\headers{Inhomogeneous Priors for Bayesian Inverse Problems}{B. Maboudi Afkham, T. Soto, M. Karamehmedovi\'c, and L. Roininen}

\title{Inhomogeneous Priors for Bayesian Inverse Problems\thanks{Submitted to the editors DATE.
\funding{The work was supported by the Research Council of Finland (RCF) through the Flagship of Advanced Mathematics for Sensing, Imaging and Modelling, and the Centre of Excellence of Inverse Modelling and Imaging (decision numbers 359183 and 353095). Babak Maboudi Afkham is also partially supported by RCF under grant numbers 371523 and Mirza Karamehmedovi\'c is supported by Villum Foundation grant number 58857.} }}

\author{Babak Maboudi Afkham\thanks{Research Unit of Mathematical Sciences, University of Oulu,
Pentti Kaiteran katu 1, Linnanmaa, Finland 
  (\email{babak.maboudi@oulu.fi}).}
\and Tom\'as Soto\thanks{School of Engineering Sciences, LUT University, Finland
  (\email{Tomas.Soto@lut.fi}, \email{Lassi.Roininen@lut.fi}).}
\and Mirza Karamehmedovi\'c\thanks{Department of Applied Mathematics and Computer Science, Technical University of Den-
mark, Kgs. Lyngby, Denmark (mika@dtu.dk).}
\and Lassi Roininen\footnotemark[3]}

\usepackage{amsopn}
